\documentclass[oneside,11pt]{article}
\usepackage{epsfig, amsfonts, amsmath, amsthm,amsbsy}
\usepackage{color}
\usepackage{subfigure}
\usepackage{hyperref}

\newtheorem{theorem}{Theorem}
\newtheorem{lemma}{Lemma}

\newtheorem{remark}{Remark}
\newtheorem{example}{Example}
\newtheorem*{example-non}{Example}

\def\ve{\varepsilon}
\def\vr{\varepsilon}
\def\ds{\displaystyle}

 \DeclareMathOperator\erfc{erfc}

 \def\cbm{\color{magenta}}

\begin{document}

 \title{A singularly perturbed convection-diffusion   parabolic problem with incompatible  boundary/initial data\thanks{This research was partially supported  by the Institute of Mathematics and Applications (IUMA), the project  PID2019-105979GB-I00  and the Diputaci\'on General de Arag\'on (E24-17R).}}
\author{J.L. Gracia\thanks{Department of Applied Mathematics, University of
Zaragoza, Spain.\  email: jlgracia@unizar.es} \and E.\ O'Riordan\thanks{School of Mathematical Sciences, Dublin City
University, Dublin 9, Ireland.\ email: eugene.oriordan@dcu.ie}}
\date{\today}
\maketitle
\begin{abstract}
A singularly perturbed parabolic problem of convection-diffusion type with incompatible inflow boundary and initial conditions is examined.
In the case of constant coefficients,  a  set of singular functions are  identified which match certain incompatibilities in the data and also satisfy the associated homogenous differential equation. When the convective coefficient only depends on the time variable and the initial/boundary data is discontinuous, then a mixed analytical/numerical approach is taken. In the case of variable coefficients and the zero level    of compatibility being  satisfied  (i.e. continuous boundary/initial data), a  numerical method is constructed whose order of convergence is shown to depend on  the next level of compatibility being  satisfied by the data.
 Numerical results are presented to support the  theoretical error bounds established for both of the approaches examined in the paper.
\end{abstract}

\section{Introduction} \label{sec:1}
Consider the following singularly perturbed parabolic problem: Find $u$ such that
\begin{subequations}\label{cont-problem-general}
\begin{align}
&Lu:=-\ve u_{xx} +a(x,t) u_x+u_t =f(x,t), \quad (x,t) \in G:= (0,1) \times (0,T], \label{GA}\\
&u(0,t) = g_L(t), \ u(1,t) = g_R(t), \ t >0,
\quad u(x,0)= \phi (x), \quad 0 \leq x \leq 1;   \label{GC}\\
&a(x,t) \geq \alpha >0, \ (x,t) \in \bar G; \quad 0 < \ve \leq 1.
\end{align}
In the case of sufficiently smooth and compatible data, the solution of this problem will contain a boundary layer of width $O(\ve)$ near the outflow boundary $x=1$. Below, we examine the issues that arise when the problem data is not sufficiently compatible at the inflow  point $(0,0)$.
If the  boundary and initial conditions are incompatible ($\phi (0) \neq g_L(0)$), a  strong \cite{B5} interior layer  will appear for small values of the singular perturbation parameter. If $\phi (0) = g_L(0)$, but the  data are still not sufficiently  compatible at $(0,0)$ then a weak \cite{B5} interior layer  appears in the solution. The path of any interior layer is located along the characteristic curve  $x=d(t)$, where $d(t)$ is implicitly defined by
\[
d'(t) = a(d(t),t), \quad  d(0)=0.
\]
To avoid the interior layer interacting with the outflow boundary, we  assume that at the final time
\begin{equation}\label{1D}
d(T) < 1.
\end{equation}
See Remark \ref{bdy-int-meet} for necessary modifications when (\ref{1D}) is not satisfied.
\end{subequations}

Let us recall the constraints on the data, $a, \, f, \, \phi, \, g_L$ and $g_R$,  for the solution $u$ to be sufficiently regular so that  classical numerical analysis is applicable; i.e.,  for $u \in  C^{4+\gamma}(\bar G)$\footnote{ As in \cite{friedman}, we define the space ${\mathcal C}^{0+\gamma}(D )$, where $D \subset \mathbf{R}^2$ is an open set, as the set of all functions that are H\"{o}lder continuous of degree $\gamma \in (0,1) $ with respect to the metric $\Vert \cdot \Vert, $
where for all ${\bf p}_i=(x_i,t_i),  \in \mathbf{R}^2, i=1,2; \
\Vert {\bf p}_1- {\bf p}_2 \Vert^2  = (x_1-x_2)^2 + \vert t_1 -t_2 \vert$.
For $f$ to be in ${\mathcal C}^{0+\gamma}(D ) $ the following semi-norm needs to be finite
\[
\lceil f \rceil _{0+\gamma , D} = \sup _{{\bf p}_1
 \neq {\bf p}_2, \ {\bf p}_1, {\bf p}_2\in D}
\frac{\vert f({\bf p}_1) - f({\bf p}_2) \vert}{\Vert {\bf p}_1- {\bf p}_2 \Vert ^\gamma} .
\] The space ${\mathcal C}^{n+ \gamma}(D ) $ is defined by
\[
{\mathcal C}^{n+\gamma }( D ) = \left \{ z : \frac{\partial ^{i+j} z}{
\partial x^i
\partial t^j } \in {\mathcal C}^{0+\gamma }(D), \ 0 \leq i+2j \leq n \right \},
\]
and $\Vert \cdot \Vert _{n + \gamma}, \ \lceil \cdot \rceil _{n+\gamma}$ are the associated norms and semi-norms.
}.
 From \cite{ladyz} we have the following result: If $a,f \in   C^{0+\gamma}(\bar G), \,  \phi \in C^{2+\gamma}[0,1], \, g_L,g_R \in C^{1+\gamma/2}[0,T]$ and
\begin{subequations}
\begin{align}
&A_0=0 \text{ with }  A_0:=g_L(0)-\phi (0), \label{C0} \\
&A_1=0 \text{ with } A_1:= -\ve \phi ''(0) +a(0,0) \phi '(0) +g_L'(0)- f(0,0), \label{C1L} \\
&g_R(0) =\phi (1); \quad  -\ve \phi ''(1) +a(1,0) \phi '(1) +g_R'(0)= f(1,0), \label{C1R}
\end{align}
then the solution of problem  (\ref{cont-problem-general}) satisfies $u \in   C^{2+\gamma}(\bar G)$. By differentiating with respect to the time variable the differential equation (\ref{GA})  and applying the above conditions on the function $u_t(x,t)$, we arrive at the following result:  If $a,f \in   C^{2+\gamma}(\bar G), \,  \phi \in C^{4+\gamma}[0,1], \, g_l,g_R \in C^{2+\gamma/2}[0,T]$ and in addition to the constraints (\ref{C0}), (\ref{C1L}) and (\ref{C1R}) we have
\begin{align}
& A_2=0 \text{ with } \nonumber \\
& A_2:=  -\ve ^2 \phi ^{(iv)}(0)+2\ve a(0,0) \phi ''' (0)- a^2 (0,0)\phi'' (0) +g_L''(0^- ) \nonumber \\
&\hspace{0.3cm} +  \ve  ( a_{xx}(0,0)  \phi'(0) +2a_x(0,0) \phi''(0)) +(a_t- aa_x)(0,0)\phi '(0) \nonumber \\
&\hspace{0.3cm} -(f_t+\vr f_{xx}-af_x)(0,0),\label{C2L} \\
&  -\ve ^2 \phi ^{(iv)}(1)+2\ve a (1,0) \phi '''(1) - a^2 (1,0) \phi''(1)  +g_R''(0^- )\nonumber \\
&\hspace{0.3cm} +  \ve( a_{xx}(1,0)    \phi'(1) +2a_x(1,0)   \phi''(1))  +(a_t- aa_x) (1,0)  \phi '(1) \nonumber \\
&\hspace{0.3cm} =(f_t+\vr f_{xx}-af_x)(1,0), \label{C2R}
\end{align}
\end{subequations}
then the solution of problem  (\ref{cont-problem-general}) satisfies $u \in   C^{4 +\gamma}(\bar G)$.

In the case of constant coefficients and when (\ref{C0}) is not satisfied, the discontinuous analytic solution in the quarter plane $x, t >0$ is given,  for example, in \cite{asymptotic3} and
in \cite{asymptotic1} an asymptotic expansion is given in the domain $\bar G$. In the case of a variable coefficient $a(t)$ and when (\ref{C0}) is  satisfied, a uniformly valid asymptotic expansion to the continuous solution of the problem posed on the quarter plane is presented in \cite{asymptotic2}.

In this paper, we examine the convection-diffusion  problem (\ref{cont-problem-general})  where the compatibility conditions (\ref{C0}), (\ref{C1L}) and  (\ref{C2L})  at $(0,0)$ are not all imposed. To  avoid additional regularity issues with the data we will assume that
$a,f \in   C^{5+\gamma}(\bar G), \, \phi \in C^7[0,1], \, g_L \in C^5[0,T], g_R \in C^3[0,T]$
and that the compatibility conditions  (\ref{C1R})  and (\ref{C2R})  at $(1,0)$ are all satisfied. We examine the problem where the initial and  left  boundary condition do not match, i.e., (\ref{C0}) is not satisfied. In this case, we
only  examine problem   (\ref{cont-problem-general}) when the convection coefficient $a(x,t)$ depends solely on the time variable. Moreover, we first separate off a singular function that matches the incompatibility at the point $(0,0)$ and then use a numerical method to approximate the difference between the solution $u$ and this singular function.

We  also examine the problem where the initial and boundary condition match, so that the zero level compatibility (\ref{C0}) is satisfied; but the higher compatibility conditions  (\ref{C1L})  and (\ref{C2L})  are not satisfied. As the solution is continuous, a numerical method can be applied directly to the problem.
If (\ref{C0}) is satisfied but the first level of compatibility  (\ref{C1L}) is not satisfied, then the order of convergence of the standard numerical method constructed in \S \ref{sec:4}  is shown to be $0.5$. If the first level of compatibility   (\ref{C0})  and (\ref{C1L}) is satisfied, then  that numerical method is essentially first order.

In \cite{cd-disc-initialA} and \cite{cd-disc-initialB},  the parabolic problem (\ref{cont-problem-general}) with compatible boundary/initial data is examined, but with a discontinuity in the initial condition $\phi(x)$ at
some internal point $x=d, \, 0 <d <1$. In this case, the  interior layer function
\[
 0.5 \erfc \left( \frac{d(t)-x}{2 \sqrt{\vr t}} \right), \ 0<d(0)< 1, \quad \erfc(z) := \frac{2}{\sqrt{\pi}} \int _{r=z} ^\infty e ^{-r^2} dr,
\]
  captures the nature of the singularity.
In \S \ref{sec:2} a set of  related functions $S_n(x,t), n \geq 0$ are constructed to model the nature of any singularity in the solution related to a lack of compatibility  between the initial and boundary condition at the point $(0,0)$. Parameter-explicit pointwise bounds on the partial derivatives of these functions are also established in \S \ref{sec:2}.  In \S \ref{sec:3}, the solution $u$ of  (\ref{cont-problem-general}) is expanded in terms of these special functions $S_n(x,t)$ as follows:
\[
u(x,t) = \sum _{i=0}^1 A_i S_i(x,t) +  \sum _{i=2}^3 B_i S_i(x,t) +v(x,t)+w(x,t),
\]
where  the amplitudes $A_i, \, i=0,1; B_i,\   i=2, 3 $ are suitably chosen so that $v,w \in   C^{4 +\gamma}(\bar G)$, where $v$ is the regular component and $w$ is the boundary layer component of the solution $u$.
With the aid of this  expansion, a numerical method is constructed in \S \ref{sec:4} to generate a numerical approximation to $u-A_0S_0$  (including also the case $A_0=0$).
The order of convergence of this method depends on whether $A_1$ is zero or not. In \S \ref{sec:5}, numerical results are presented for  sample test  problems to illustrate the performance of the method and to validate the orders of convergence established in the two main Theorems \ref{th_a(t)} and \ref{th_a(x,t)} in \S \ref{sec:4}.  Technical details associated with establishing bounds on the derivatives of the functions $S_n(x,t)$ are given in the appendix.

{\bf Notation:} Throughout the paper, $C$ denotes a generic constant that is independent of the singular perturbation parameter $\vr$ and all the discretization parameters. The $L_\infty$ norm on the domain $D$ will be denoted by $\Vert \cdot \Vert_D$. We also define
 the following interior layer function
\[
 E_\gamma(x,t):= e^{-\frac{\gamma(x-d(t))^2}{4\ve t }}, \quad 0 <\gamma \le 1.
\]
 If $\gamma=1$, we simply write $E_1(x,t)=E(x,t)$.

\section{A set of singular functions with incompatiblities} \label{sec:2}

We  now  define a set of singular functions  $S_n(x,t)$, which will form  a basis for the regularity expansion of the solution  $u(x,t)$ of problem \eqref{cont-problem-general}, which is constructed  in  Theorem \ref{theorem:decomp} (for $A_0\ne 0$ and $a=a(t)$) and Theorem \ref{theorem:decomp-general} (for $A_0= 0$ and $a=a(x,t)$).
For all $n \geq 0$:
\begin{equation}\label{Sn-def}
S_n(x,t) := \frac{\psi _n^+(x,t) +(-1)^n \psi ^-_n(x,t)}{a^n(0,0)},
\end{equation}
where the functions $\psi ^\pm_n(x,t), n \geq -1$ are defined by
\begin{subequations}\label{def-basic}
\begin{align}
\psi ^-_n (x,t) &:= (-1)^n 2^{n-1} n! (\ve t )^{n/2} \erfc _n (\chi ^- (x,t)), \\
\psi ^+ _n (x,t) &:= (-1)^n 2^{n-1} n! (\ve t )^{n/2} e^{\frac{xd(t)}{\ve t}} \erfc _n (\chi ^+ (x,t)), \\
 \chi ^\pm (x,t)& := \frac{x\pm d(t)}{2\sqrt{\ve t }}, \quad \psi _{-1}^\pm (x,t) := -\frac{E(x,t)}{2\sqrt{\ve \pi t}},
\end{align}
\end{subequations}
and the iterated complementary error functions are
\[
\erfc_{-1}(x):=\frac{2}{\sqrt{\pi}}e^{-x^2}, \quad \erfc_n(x):=\int_{s=x}^\infty\erfc_{n-1}(s) \, ds, \quad n \geq 0.
\]
Observe that the first function $S_0$ is discontinuous and
\[
S_1 \in C^{0+\gamma} (\bar G), \qquad S_{2n},S_{2n+1}  \in C^{2n+\gamma} (\bar G), \quad  n \geq 1.
\]
 In the next lemma, we establish bounds on the derivatives of the first three functions $S_n, \, n=0,1,2.$ These bounds indicate both the strength of the singularity at $t=0$ and how certain derivatives can depend on inverse powers of $\vr$.

\begin{lemma} \label{lemma:BoundsSi}
The  function $S_0(x,t)$ satisfies the bounds
\begin{subequations}\label{S0-bounds}
\begin{align}
&  \vert S_0 \vert \leq C, \quad   \left \vert \frac{\partial^i S_0}{\partial t^i} \right \vert   \leq C \left[ \frac{1}{t} \left(1+\sqrt{\frac{t}{\ve }} \right)\right ]^i E_\gamma(x,t), \quad i=1,2,
\\
& \left \vert \frac{\partial^i S_0}{\partial x^i } \right \vert   \leq \frac{C} {\vr^i} \left(\frac{\vr}{t}+ \left(\frac{\vr}{t} \right)^{i/2} \right) E_\gamma(x,t), \quad i=1,2,3;
\end{align}
\end{subequations}
the  function $S_1(x,t)$ satisfies
\begin{subequations}\label{S1-bounds}
\begin{align}
& \vert S_1 \vert \leq  C, \quad
\left \vert \frac{\partial S_1}{\partial t  } \right \vert \leq  C, \ \left \vert \frac{\partial ^2S_1}{\partial t ^2} \right \vert  \leq C\frac{1}{t} \left(1+\sqrt{\frac{t}{\ve }} \right) E_\gamma(x,t)  +C,
\\
& \left \vert \frac{\partial S_1}{\partial x } \right \vert \leq   C, \  \left \vert \frac{\partial ^2 S_1}{\partial x ^2} \right \vert \leq  \frac{C}{\ve} E_\gamma(x,t)  +C, \  \left \vert \frac{\partial ^3 S_1}{\partial x ^3} \right \vert \leq  \frac{C}{\vr \sqrt{\vr t}} E_\gamma(x,t)   +C;
\end{align}
\end{subequations}
and the function $S_2(x,t)$ satisfies
\begin{subequations}\label{S2-bounds}
\begin{align}
 & \vert S_2 \vert \leq   C, \quad  \left \vert \frac{\partial S_2}{\partial x}\right \vert   \leq    C,
 \\
& \left \vert \frac{\partial S_2}{\partial t  } \right \vert \leq    C,
\
 \left \vert \frac{\partial ^2 S_2}{\partial t ^2} \right \vert
\leq  C \left(1+\frac{\vr}{t} \right)  \left (1+\sqrt{\frac{t}{\vr}} \right) E_\gamma(x,t) +C,
\\
& \left \vert \frac{\partial ^2 S_2}{\partial x ^2} \right \vert \leq   C \left(1+\sqrt{\frac{t}{\ve }} \right) E_\gamma(x,t),
\\
& \left \vert \frac{\partial ^3 S_2}{\partial x ^3} \right \vert \leq   \frac{C}{\ve}
\left(1+\sqrt{\frac{t}{\ve}}+\sqrt{\frac{\ve}{t}}\right)E_\gamma(x,t).
\end{align}
\end{subequations}
\end{lemma}
\begin{proof} In the appendix, bounds on the partial derivatives of the functions $\psi ^\pm_n(x,t)$ are established. These are used to  prove the bounds on $S_n$. The bounds on $S_0$  follow directly from (\ref{bounds-singular-functions}) and (\ref{bounds-star-functions}).
Using (\ref{initial-recurrance}) and the recurrence relation  (\ref{deriv-relations-general}) we deduce the following
\begin{align*}
\frac{\partial S_1}{\partial t}  &= \frac{a(d(t),t)}{a(0,0)} S_0 + x \frac{p(t) \psi ^+_1}{\ve t^2} , \quad  p(t):= ta(d(t),t)-d(t);
\\
a^2(0,0) \frac{\partial S_2}{\partial t}  &=2(\ve S_0 - a(0,0)a(d(t),t) S_1 + \frac{d(t)}{t} \psi ^+_1)
 - \frac{p(t)}{\ve t^2} (d(t) \psi ^+_2 - \psi ^+_3).
\end{align*}
The bounds on the time derivatives of $S_1$ and $S_2$  follow. To deduce the bounds on the space derivatives of these components,  we first note that from \eqref{deriv-relations-general}, we have
\begin{align*}
a(0,0)\frac{\partial S_1}{\partial x} & = \psi ^+_0 - \psi ^- _0 + \frac{ d(t)}{\ve t} \psi _1^+,
\\
a(0,0) \frac{\partial ^2S_1}{\partial x ^2} & = \frac{d(t)}{t \ve}\left( 2 \psi ^+_0 +  \frac{d(t)}{t \ve} \psi ^+_1 \right),
 \\
a(0,0) \frac{\partial ^3S_1}{\partial x ^3} &  = \frac{d(t)}{t \ve}\left( 2 \frac{\partial  \psi ^+_0}{\partial x}  +  \frac{d(t)}{t \ve} \left(\psi ^+_0 +\frac{d(t)}{t \ve} \psi ^+_1\right) \right),
\end{align*}
and from  \eqref{bounds-psi12} the bounds on the space derivatives of $S_1$ follow.
Next, we deduce bounds on the space derivatives of $S_2$.
From the definitions in   \eqref{def-basic}, the following bounds are obtained
\[
  \left \vert \frac{\partial ^i \psi ^-_j}{\partial x ^i} \right \vert \le C (\sqrt{\ve t})^{j-i} E_\gamma(x,t)     + C,\quad  j=1,2; \quad i=1,2,3,
\]
and using  the recurrence relation  \eqref{deriv-relations-general} we get the bounds
\begin{subequations}
\begin{align}
& \left \vert \frac{\partial \psi ^+_1}{\partial x}\right \vert \leq C E_\gamma(x,t), \
\left \vert \frac{\partial^2 \psi ^+_1}{\partial x^2}\right \vert \leq \frac{C}{\vr}\left(1+\sqrt{\frac{\vr}{t}} \right)E_\gamma (x,t), \label{psi1-derivatives}
\\
 & \left \vert \frac{\partial^3 \psi ^+_1}{\partial x^3}\right \vert \leq \frac{C}{\vr t}\left(1+\sqrt{\frac{t}{\vr}} \right)E_\gamma (x,t),
\\
& \left \vert \frac{\partial \psi ^+_2}{\partial x}\right \vert \leq C \sqrt{\ve t} E_\gamma(x,t), \quad  \left \vert \frac{\partial ^2 \psi ^+_2}{\partial x^2}\right \vert \leq C \left(1+\sqrt{\frac{t}{\vr}} \right) E_\gamma(x,t),  \\
& \left \vert \frac{\partial ^3 \psi ^+_2}{\partial x^3}\right \vert \leq  \frac{C}{\ve}\left(1+\sqrt{\frac{t}{\vr}}  +\sqrt{\frac{\vr}{t}} \right) E_\gamma(x,t).
\end{align}
\end{subequations}
 The bounds on the space derivatives of $S_2$ follow immediately from the  bounds above on the space derivatives of the singular functions $\psi ^-_i$ and $\psi_i^+$.

\end{proof}

Observe that the strength of the singularity at $(0,0)$ in each of the functions $S_n$ weakens as $n$ increases.
Using the identities in (\ref{deriv-relations-general}), we can  deduce bounds on the remaining functions $S_n, \, n \ge 3$:
\begin{subequations}\label{Sn-bounds}
\begin{align}
& \left \vert  \frac{\partial  S_{n}}{\partial x }  \right \vert \leq C ,\ n \geq 2, \quad \left \vert  \frac{\partial  ^2S_{n}}{\partial x ^2}  \right \vert \leq C ,\ n \geq 4, \quad \left \vert  \frac{\partial ^3 S_{n}}{\partial x ^3}  \right \vert \leq C ,\ n \geq 6,
\\
& \left \vert  \frac{\partial ^2 S_{3}}{\partial x ^2}  \right \vert \leq C\left(1+t\left(1+\sqrt{\frac{t}{\ve}}\right) E_\gamma(x,t) \right) ,
\\
& \left \vert  \frac{\partial ^3 S_{3+n}}{\partial x ^3}  \right \vert \leq C\left(1+t^n\left(1+\sqrt{\frac{t}{\ve}}\right)^{3-n}E_\gamma(x,t)\right) ,\quad n=0,1,2,
\\
& \left \vert  \frac{\partial  S_{n}}{\partial t }  \right \vert \leq C ,\ n\geq 2, \
\left \vert  \frac{\partial ^2 S_{3}}{\partial t ^2}  \right \vert \leq C \left(1+\sqrt{\frac{\ve}{t}} E_\gamma(x,t)\right),
 \left \vert  \frac{\partial ^2 S_{n}}{\partial t ^2}  \right \vert \leq C , \  n\geq 4.
\end{align}
\end{subequations}
\section{The continuous problem} \label{sec:3}

In the following result the asymptotic behaviour of  the solution $u$ to problem~\eqref{cont-problem-general} is given when the convective coefficient $a$ depends only on the time variable and $u$ is discontinuous at $(0,0)$.

\begin{theorem} \label{theorem:decomp}
 Assume that $a(x,t)=a(t), \, \forall (x,t) \in \bar G$ and  $a_t(0)=0$. The solution $u$ of (\ref{cont-problem-general}) can be expanded as follows
\begin{subequations}
\begin{equation}\label{decomp}
u(x,t) = \sum _{i=0}^1 A_i S_i(x,t)  +\sum _{i=2}^3 B_i S_i(x,t)  +v(x,t)  + w(x,t),
\end{equation}
where the constants $A_i$ are defined in \eqref{C0} and \eqref{C1L}. The constants   $B_2,B_3$ are defined such that  for $0\le i+2j\le 4$
\begin{align}
\left \vert \frac{\partial ^{i}v}{\partial x^i} \right \vert & \leq C(1+ \ve^{2-i}), \quad \left\vert \frac{\partial ^{j}v}{\partial t^j} \right \vert \leq  C ; \label{bnds-v}
\\
\left\vert \frac{\partial ^{i+j}w}{\partial x^i\partial t^j} \right \vert & \leq C \ve^{-i}(1+\ve^{1-j})e^{-\alpha (1-x)/\ve}. \label{bnds-w}
\end{align}
\end{subequations}
\end{theorem}
\begin{proof}
With the assumptions $a=a(t), \, a_t(0)=0$ and noting  (\ref{def-p}), we have
\begin{equation}\label{assume-a(t)}
 LS_i = p(t)\frac{ \psi ^+_{i+1}}{\ve t^2},\  p(t)=t^3  P(t), \  \vert P(t) \vert \leq C \quad \hbox{and}\quad LS_0 \in   C^{2+\gamma}(\bar G).
\end{equation}
We identify  the remainder $R$ by
\[
 R:= u(x,t) - \sum _{i=0}^1 A_i S_i(x,t), \quad  R(0,t) = g_L(t)-C_L(t), \ R(x,0) =\phi (x),
\]
where
\[
 C_L(t):=  \sum _{i=0}^1 A_i S_i(0,t)= \sum _{i=0}^1 A_i \left( \frac{d(t)}{a(0)}\right)^i.
\]
Then $C_L(0)=A_0, \ C_L'(0)=A_1$ and the remainder function $R$ satisfies
\begin{equation} \label{LRemainder}
 LR= f-\frac{ p(t)}{\ve t^2} \left(A_0\psi ^+_1 +\frac{A_1}{a(0)} \psi ^+_2\right) \in   C^{2+\gamma}(\bar G).
\end{equation}
 Hence, $R \in   C^{2 +\gamma}(\bar G)$ as the amplitudes $A_i,i=0,1$ have been chosen so that the compatibility conditions (\ref{C0}) and (\ref{C1L})  are satisfied by the problem data defining $R$.

The remainder $R$ is further decomposed as follows
\[
R(x,t)= \sum _{n=2}^5 B_n S_n(x,t)+z(x,t)+v_S(x,t)+w(x,t),
\]
 with $z, v_S, w\in C^{4+\gamma}(\bar G)$.  The regular component $v$ of the remainder $R$ will be determined by $v=\sum _{n=4}^5 B_n S_n(x,t)+z+v_S$ and the functions $z$ and $v_S$ are required in our decomposition due to the weak singular right-hand side of the differential equation~\eqref{LRemainder}. The boundary layer function $w$ of $R$ will satisfy  the problem
$Lw=0, w(0,t)=w(x,0)=0, w(1,t) \neq 0$. All these functions and the constants $B_n$ are specified below.


Consider the following function
\[
z(x,t) : =  \phi(x)  +z_0(x,t)+\ve z_1(x,t) +\ve ^2 R_z(x,t);
\]
where
\begin{align*}
& L_0z_0 =f+\ve \phi ''(x) - a(t) \phi'(x), \ 0 <x \le 1, \ t>0,
 \\
 &  z_0(0,t)  = g_L(t) -C_L(t) - \phi (0)- \sum _{n=2}^5 B_n S_n(0,t), \ t>0, \  z_0(x,0)=0, \ 0 \le x \le 1,
 \\
& L_0 z_1  = \frac{\partial ^2 z_0}{\partial x ^2}, \ 0 <x \le 1, \ t>0, \  z_1(0,t) =0, \ t>0, \ z_1(x,0)=0, \   0 \le x \le 1,
\\
& LR_z  = \frac{\partial ^2 z_1}{\partial x ^2},  \ (x,t)\in G, \ R_z(0,t) =R_z(1,t) =0, \ t>0, \ R_z(x,0)=0, \ 0 \le x \le 1,
\end{align*}
with the reduced differential operator
\[
 L_0 := \frac{\partial}{\partial t} + a(t) \frac{\partial}{\partial x}.
\]
From this construction, $z(1,t) =z_0(1,t)+\ve z_1(1,t)$ and
\begin{eqnarray*}
Lz=f, \quad z(0,t)= R(0,t) - \sum _{n=2}^5 B_n S_n(0,t), \ z(x,0)=R(x,0).
\end{eqnarray*}
Note that the function $z_0$ satisfies $z_0(0^+,0)=z_0(0,0^+)$ and the first level compatibility condition
\[
f(0,0) + \ve \phi ''(0) - a(0) \phi '(0)= g_L'(0) - C_L'(0),
\]
is satisfied automatically. Hence,  from \cite{bobisud} and \cite{linss+stynes},   the function $z_0$ belongs to the space{\footnote{The space ${\mathcal C}^{n,\gamma}(D ) $ is defined by
\[
{\mathcal C}^{n,\gamma }( D ) = \left \{ z : \frac{\partial ^{i+j} z}{
\partial x^i
\partial t^j } \in {\mathcal C}^{0+\gamma }(D), \ 0 \leq i+j \leq n \right \},
\]}}
$C^{1,\gamma} (\bar G),$  due to the particular   definition of $z_0$.

 Now the parameters $B_n, n=2,3,4,5$ are chosen so that the necessary compatibility conditions on the reduced solution $z_0$ are imposed in order  that $z_0 \in C^{5,\gamma} (\bar G)$. Then,
$z_1 \in C^{4,\gamma} (\bar G)$ and $R_z\in  C^{4+\gamma} (\bar G)$. Hence, $z\in  C^{4+\gamma} (\bar G)$ and
\[
\left \vert \frac{\partial^{i+j}z}{\partial x^i\partial t^j}\right \vert \le C \left(1+\vr^{2-i}\right), \quad 0\le i+2j\le 4.
\]

We next define the component $v_S$ as the solution of the initial-boundary value problem
\begin{align*}
 L v_S & = LR-Lz-\sum _{n=2}^5B_nLS_n= -\sum _{n=0}^1A_nLS_n-\sum _{n=2}^5B_nLS_n, \quad (x,t) \in G, \\
v_S(0,t) & = v_S(1,t)=0, \, t>0, \quad  v_S(x,0) = 0, \, 0 \le x \le 1.
\end{align*}
Observe that  $v_S \in C^{4+\gamma} (\bar G)$ as
\[
(Lv_S)(0,0)=(Lv_S)_x(0,0)=(Lv_S)_{xx}(0,0)=(Lv_S)_t(0,0)=0,
\]
and $ Lv_S \in C^{2+\gamma} (\bar G)$.
To deduce bounds on the derivatives of the component $v_S$ consider the stretched variables
\begin{equation} \label{stretching}
\tau = \frac{t}{\ve}, \ \zeta = \frac{x}{\ve}
\end{equation}
and  we denote $\tilde g (\zeta, \tau ) := g(x,t)$ for any function $g$. Then, we have
\[
\frac{x+d(t)}{2\sqrt{\ve t}} = \frac{\zeta+\tilde d_1(\tau)}{2\sqrt{\tau}}, \ \text{ with } \ \tilde d_1(\tau) := \frac{1}{\ve} \int _{s=0}^{\ve \tau} a(s) \ ds = \int _{s=0}^{ \tau}  \tilde a(s) \ ds.
\]
Using (\ref{assume-a(t)}),  we have that
\[
- \frac{\partial ^2 \tilde v_S}{\partial \zeta ^2} + \tilde a(\tau) \frac{\partial  \tilde v_S}{\partial \zeta } + \frac{\partial  \tilde v_S}{\partial \tau } =  - \ve^2  \tilde \Phi ^+(\zeta , \tau),
\]
where
\[
 \tilde \Phi ^+(\zeta , \tau) :=  \tau  \tilde P(\tau) \left(\sum _{n=0}^1 \frac{A_n \tilde \psi ^+_{n+1}}{\vr} +\sum _{n=2}^5  \frac{B_n \tilde \psi ^+_{n+1}}{\vr} \right).
\]
Then, from the definition (\ref{def-basic}) of the basic functions $ \psi^- _n$ and  $\psi ^+_n$, we have
\[
\left \vert \frac{\partial ^{i+j} \tilde \Phi ^+}{\partial \zeta  ^i\partial  \tau ^j} \right \vert \leq C, \quad 0 \leq i+2j \leq 2.
\]
From \cite{friedman} and \cite{ladyz}, we have the following estimates for the partial derivatives of $\tilde v_S$
\[
\left\vert \frac{\partial ^{i+j} \tilde v_S}{\partial \zeta  ^i\partial  \tau ^j} \right \vert \leq C \ve ^2, \quad 0 \leq i+2j \leq 4.
\]
Returning to the original variables, we get that
\[
\left\vert \frac{\partial ^{i+j}  v_S}{\partial x ^i\partial t ^j} \right \vert \leq C \left(1+ \ve ^{2-(i+j)}\right), \quad 0 \leq i+2j \leq 4.
\]
We can now define the regular component $ v \in C^{4+\gamma}(\bar G)$ to be
\[
v(x,t):= z(x,t) +v_S(x,t) +\sum _{n=4}^5B_n S_n(x,t).
\]
Finally, consider the boundary layer component $w$; it is the solution of
\begin{align*}
L w & =0, \quad (x,t) \in G, \\
w(0,t) &=0, \ w(1,t)=\left(R-v-\sum _{n=2}^3B_nS_n\right)(1,t),\ t>0, \\
w(x,0) &=0, \ 0 \le x \le 1.
\end{align*}
The bounds on $w$ are established as in \cite[Theorem 1]{math-comp-2016}. We note that we require the assumption \eqref{1D}  to establish the bounds on the derivatives of the boundary layer function  $w$. This assumption guarantees that $d(T) <1$ and then the interior and boundary layers do not interact  with each other.

\end{proof}
In the next theorem, we consider the case where (\ref{C0}) is satisfied and the solution is continuous. In this case we can relax the constraints on the coefficient $a(x,t)$  and allow this coefficient to vary in  both space and time.

\begin{theorem} \label{theorem:decomp-general}
 Assume that $a_x(0,0)=0$ and $g_L(0)=\phi (0)$. The solution $u$ of (\ref{cont-problem-general}) can be expanded as follows:
\begin{subequations}
\begin{equation}\label{decomp-gen}
u(x,t) = A_1S_1(x,t)+\sum _{i=2}^3 B_i S_i(x,t)  +v(x,t) +w(x,t),
\end{equation}
where $A_1$ is defined in \eqref{C1L} and for $0 \leq i+2j \leq 4$
\begin{align}
\left\vert \frac{\partial ^{i}v}{\partial x^i} \right \vert &\leq  C(1+ \ve^{2-i}), \quad \left\vert \frac{\partial ^{j}v}{\partial t^j} \right \vert \leq  C, \label{bnds-v-gen}
\\
\left\vert \frac{\partial ^{i+j}w}{\partial x^i\partial t^j} \right \vert  &\leq  C \ve^{-i}(1+\ve^{1-j})e^{-\alpha (1-x)/\ve}. \label{bnds-w-gen}
\end{align}
\end{subequations}
\end{theorem}
\begin{proof} Follow the argument in Theorem \ref{theorem:decomp-general}, but now we define the remainder to be
$
R:= u(x,t) - A_1 S_1(x,t)$
which  satisfies
\begin{equation} \label{LRemainder-2}
 LR=f-\frac{p(t)}{\ve t^2} \frac{A_1}{a(0,0)} \psi ^+_2   - (a(x,t) - a(d(t),t)) A_1\frac{\partial S _1}{\partial x}.
\end{equation}
 We examine the regularity of the function $R$. Compare \eqref{LRemainder-2} with \eqref{LRemainder}. As in Theorem~\ref{theorem:decomp}, the term $f-\frac{p(t)}{\ve t^2} \frac{A_1}{a(0,0)} \psi ^+_2$ in the right-hand side of \eqref{LRemainder-2} belongs to $C^{2+\gamma}(\bar G)$ since $\psi ^+_2(x,t)  \in C^{2+\gamma}(\bar G)$. We now consider the other term $(a(x,t) - a(d(t),t)) A_1\frac{\partial S _1}{\partial x}$ of \eqref{LRemainder-2}. Observe that, if $a_x(0,0)=0$ then
\begin{align*}
 a(x,t)-a(d(t),t) & = \int _{s=d(t)}^x a_x(s,t) \, ds \\
&= \int _{s=d(t)}^x  \int _{r=0}^t  a_{xt}(s,r) \, dr \, ds +  \int _{s=d(t)}^x  \int _{r=0}^s  a_{xx}(r,0) \, dr \, ds.
\end{align*}
Hence, 
\begin{align*}
\vert  a(x,t)-a(d(t),t)\vert &\leq  Ct\vert x-d(t) \vert  + C\vert (x-d(t))\vert \, \vert  (x+d(t))\vert \\
 &\leq  Ct\vert x-d(t) \vert  +C(x-d(t))^2.
\end{align*}
 In addition, we  have that
\begin{align*}
 a(0,0) \frac{\partial}{\partial t} \left( \frac{\partial S_1}{\partial x} \right) & =  \frac{\partial}{\partial t} \left( \psi ^+_0 -\psi _0 \right)  +  \frac{\partial}{\partial t} \left( \frac{d(t)}{\ve t} \psi ^+_1  \right)  \quad \hbox{and}
\\
 \frac{\partial }{\partial t} \left( \psi ^+_0- \  \psi _0 \right) & =\frac{d(t)-2ta(d(t),t)}{2t\sqrt{\ve \pi t }} E(x,t) + \frac{p(t)x}{\ve t^2} \psi ^+_0,
 \\
\frac{\partial }{\partial t} \left( \frac{d(t)}{\ve t} \psi ^+_1 \right) & =\psi ^+_1 \frac{\partial }{\partial t} \Bigl( \frac{d(t)}{\ve t}   \Bigr) +  \frac{d(t)}{\ve t}\frac{\partial \psi ^+_1}{\partial t}.
\end{align*}
Recall~\eqref{tdpsi+dt} and observe also that
\[
(x-d(t))(x+d(t))\frac{1}{t}  \psi^+ _0  \in C^{0+\gamma}(\bar G).
\]
Together, these imply that
\[
(x-d(t))t \frac{\partial \psi ^+_1}{\partial t}, \ (x-d(t))^2 \frac{\partial \psi ^+_1}{\partial t} \in C^{0+\gamma}(\bar G).
\]
 Therefore,
\[
(x-d(t))t \frac{\partial ^2 S_1}{\partial x \partial t} , \, (x-d(t))^2 \frac{\partial ^2 S_1}{\partial x \partial t} \in C^{0+\gamma}(\bar G).
\]
 Thus, it follows that  $LR \in   C^{2+\gamma}(\bar G)$ if we assume that $a_x(0,0)=0$. Hence, $R \in   C^{2 +\gamma}(\bar G)$.
 Furthermore, using  the stretched variables $\tau$ and $\zeta$  defined in ~\eqref{stretching}, we note that if
\[
\Phi _1(x,t) := \frac1{a(0,0)} (a(x,t)-a(d(t),t))  \left( (\psi ^+_0 -  \psi^- _0) +\frac{d(t)}{\ve} \psi ^+_1 \right),
\]
then, as $\vert \tilde a(\zeta, \tau) - \tilde a(d(\tau), \tau) \vert \leq C \ve ^2(  \zeta +\tau )^2$, it follows that
\[
\vert \tilde \Phi _1(\zeta ,\tau )  \vert  \leq C \ve ^2
  (  \zeta +\tau )^2 \left( \vert \tilde \psi ^+_0 \vert + \vert \tilde  \psi^- _0 \vert + \tau \vert \tilde \psi ^+_1 \vert \right),
\]
 Use this expression and \eqref{bounds-psi12} to deduce bounds on the derivatives of the corresponding component $v_S$ of the solution $u$.
The argument is then completed as in the proof of the previous theorem.
\end{proof}

 In the next section, we describe a numerical method that will generate a numerical approximation to $y= u- A_0 S_0$. If $A_0=0$ (i.e., the zero level compatibility is satisfied), note that $y=u$. The function $y$ satisfies the singularly perturbed problem
\begin{subequations}\label{problem2}
\begin{align}
&Ly =f-A_0LS_0, \quad (x,t) \in G, \\
&y(0,t) = g_L(t)-A_0, \ y(1,t) = g_R(t)-A_0S_0(1,0), \quad t >0,\\
&y(x,0)= \phi (x), \quad 0 \leq x \leq 1.
\end{align}
\end{subequations}

\section{Numerical method} \label{sec:4}

Let  $N$ and $M=O(N)$ be two positive integers. We approximate problem~\eqref{problem2} with a finite difference scheme on a mesh  $ \bar G^{N,M}=\{x_i\}^N_{i=0} \times \{t_j\}_{j=0}^M $.  We denote by $\partial G^{N,M}:=\bar G^{N,M}\backslash G.$
The mesh $\bar G^{N,M}$ incorporates a uniform mesh  ($t_j:=k j$ with $k=T/M$) for the time variable and a piecewise-uniform mesh for the space variable with $h_i:= x_i-x_{i-1}$.  The piecewise uniform mesh $\{x_i\}^N_{i=0}$ is a Shishkin mesh \cite{fhmos} which splits the interval $[0,1]$ into the two subintervals
\[
[0,1-\sigma]\cup [1-\sigma ,1], \quad \hbox{where} \quad \sigma : =\min \left \{0.5, \frac{\vr}{\alpha}  \ln N \right\}.
\]
The $N$ space mesh points are distributed in the ratio $N/2:N/2$ across the two subintervals.
The discrete problem\footnote{We use the following notation for the finite difference approximations of the derivatives:
\begin{align*}
&D^-_t Y (x_i,t_j) := \ds\frac{Y(x_i,t_j)-Y(x_i,t_{j-1})}{k}, \quad
D^-_x Y(x_i,t_j) :=\ds\frac{Y(x_i,t_j)-Y (x_{i-1},t_j)}{h_i}, \\
& D^+_x Y (x_i,t_j)  :=\ds\frac{Y(x_{i+1},t_j)-Y(x_i,t_j)}{h_{i+1}}, \
 \delta^2_x Y(x_i,t_j) := \ds\frac{2}{h_i+h_{i+1}}(D_x^+Y(x_i,t_j)-D^-_x Y(x_i,t_j)).
 \end{align*}}
 is: Find $Y$ such that
\begin{subequations}\label{discr-probl}
\begin{align}
L^{N,M} Y &:=  -\ve \delta ^2_x  Y  + a D^-_x Y + D^-_t Y =f-A_0LS_0,\quad t_j >0,  \\
Y (x_i,0) & = y (x_i,0), \ 0 <x_i<1,\\  Y(0,t_j) &=y(0,t_j), \quad  Y(1,t_j) = y(1,t_j),\ t _j \geq 0.
\end{align}
\end{subequations}

We form a global approximation $ \bar Y$ using simple  bilinear interpolation:
\[
\bar Y(x,t) : = \sum _{i=0,j=1}^{N,M} Y(x_i,t_j) \varphi _i(x)  \eta_j(t),
\]
where $\varphi _i(x)$ is the standard hat function centered at $x=x_i$ and $\eta _j(t) :=(t-t_{j-1})/k, \, t \in (t_{j-1},t_j], \,  \eta _j(t) :=0$ otherwise.

 In the next theorem, a convergence result is given in the particular case of $a=a(t)$ and $g_L(0) \ne \phi (0)$.  In this case, the solution of problem~\eqref{cont-problem-general} is decomposed as in Theorem~\ref{theorem:decomp} and after separating off the singular function $S_0$, the numerical method  (\ref{discr-probl}) is applied to approximate $y=u-A_0S_0$.

\begin{theorem} \label{th_a(t)}  Assume that $a(x,t)=a(t), \forall (x,t) \in \bar G, a_t(0)=0$.and $ M=O(N)$.
If $Y$ is the solution of (\ref{discr-probl}) and $y$ is the solution of (\ref{problem2}), then
\begin{eqnarray*}
\Vert   \bar Y- y \Vert _{\bar G} \leq  C \vert A_1 \vert N^{-1/2}+C N^{-1} \ln N.
\end{eqnarray*}
\end{theorem}
\begin{proof}
As in the case of the continuous problem, the discrete solution can be decomposed into the sum $Y =A_1S_1^N +A_2 S_2^N+ V+W$, where
\begin{align*}
L^{N,M}V &=Lv,\ (x_i,t_j) \in {G^{N,M}} \hbox{ and }  V= v,\ ({ x_i},t_j) \in  { \partial G^{N,M};} \\
 L^{N,M}W &=0,\ (x_i,t_j) \in { G^{N,M}} \hbox{ and }   W= w,\ ({x_i},t_j) \in  { \partial G^{N,M};} \\
  L^{N,M}S_k^N&=LS_k,\ (x_i,t_j) \in { G^{N,M}}  \hbox{ and }   S_k^N= S_k,\ ({ x_i},t_j) \in  { \partial G^{N,M},} \ k=1,2.
\end{align*}
Using the bounds on the derivatives (\ref{bnds-w}) of the component $ w$  to obtain appropriate truncation error  estimates, the discrete maximum principle with  a suitable discrete barrier function and following the arguments in \cite{ria}, we can establish the following bounds
\begin{equation} \label{error-w}
 \vert ( w- W)(x_i,t_j) \vert  \leq  C N^{-1} \ln N , \quad   (x_i,t_j)\in\bar G^{N,M}.
\end{equation}
The error due to the regular component $ v$ can be bounded in a classical way \cite{mos2} to deduce that
\begin{equation} \label{error-v}
 \vert ( v- V)(x_i,t_j) \vert  \leq  C N^{-1} , \quad   (x_i,t_j)\in\bar G^{N,M}.
\end{equation}
Let us now consider the two weakly singular functions $S_1,S_2$ and their numerical approximations $S^N_1, S_2^N$.
For both functions, the truncation error is denoted by
 \[
{\cal T}_{S_k;i,j}:= L^{N,M} ( S_k -S^N_k)(x_i,t_j),
\]
then
\begin{align*}
 \vert {\cal T}_{S_k;i,j}\vert  &\leq C \ve (h_i+h_{i+1}) \left\Vert \frac{\partial ^3  S_k (x,t_j)}{\partial x^3} \right\Vert _{(x_{i-1},x_{i+1}) } \\
 &+C \min \left\{ h_i\left\Vert \frac{\partial ^2  S_k(x,t_j)}{\partial x^2} \right\Vert _{(x_{i-1},x_{i}) }, \left\Vert \frac{\partial   S_k(x,t_j)}{\partial x} \right\Vert _{(x_{i-1},x_{i}) }  \right \}
 \\
 &+
C\min \left\{ \frac{1}{k} \int _{w=t_{j-1}}^{t_j} \int _{r=w}^{t_j}  \left\vert \frac{\partial ^2  S_k(x_i,r)}{\partial t^2} \right\vert dr \ dw ,    \left\Vert \frac{\partial   S_k(x_i,t)}{\partial t} \right\Vert _{  (t_{j-1}, t_j )}  \right\},
\end{align*}
as
\[
\vert D_t^-  S_k(x_i,t_j) \vert  \leq \frac{1}{k} \int _{r=t_{j-1}} ^{t_j} \left \vert \frac{\partial  S_k(x_i,r)}{\partial r}   \ dr \right \vert \leq C \left \Vert \frac{\partial   S_k(x_i,t)}{\partial t} \right \Vert _{  (t_{j-1}, t_j )}.
\]
Note also that at each time level,
\[
\left(-\ve \delta ^2_x    + a D^-_x  +\frac{1}{k}I \right)  ( S_k -S_k^N)(x_i,t_j) =  {\cal T}_{S_k;i,j} + \frac{1}{k}  (S_k -S_k^N)(x_i,t_{j-1}), \ t_j >0 .
\]
In the case of the weaker  singular function $S_2$, we use the bounds (\ref{S2-bounds})  so that  the truncation error at the first time level $t=t_1$ is
\begin{align*}
\vert {\cal T}_{S_2;i,1} \vert & \le  \vr \left \Vert  \frac{\partial^2 S_2}{\partial x^2} \right \Vert_{(x_{i-1},x_{i+1})}
+ a \left \Vert  \frac{\partial S_2}{\partial x} \right \Vert_{(x_{i-1},x_i)}
+  \left \Vert  \frac{\partial S_2}{\partial t} \right \Vert_{(t_0,t_1)} \\
& \le C (\ve +\sqrt{\ve t_1})  e^{-\gamma \frac{(x_i-at_1)^2}{4\vr T}} \le C.
\end{align*}
At the next time levels $t_n, \, n\ge 2$, we  again use the bounds (\ref{S2-bounds}) to deduce the truncation error bounds
\begin{align*}
& \vert {\cal T}_{S_2;i,j} \vert  \le C \vr (h_i+h_{i+1}) \left \Vert  \frac{\partial^3 S_2}{\partial x^3} \right \Vert_{(x_{i-1},x_{i+1})}
+a h_i \left \Vert  \frac{\partial^2 S_2}{\partial x^2} \right \Vert_{(x_{i-1},x_i)}
+ C k \left \Vert  \frac{\partial^2 S_2}{\partial t^2} \right \Vert_{(t_{j-1},t_j)} \\
& \le C N^{-1} \left( 1+\sqrt{\frac{\vr}{t_j}} + \sqrt{\frac{t_j}{\vr}} \right)E_\gamma(x,t)+ CM^{-1} \left(1+\sqrt{\frac{\vr}{t_{j-1}}} +\frac{\vr}{t_{j-1}}  + \sqrt{\frac{t_j}{\vr}} \right)   E_\gamma(x,t)\\
&  \le C  \left(\frac{M^{-1}}{\sqrt{\ve}}+\frac{M^{-1/2}}{\sqrt{j-1}}+\frac{\vr}{j-1} \right) E_\gamma(x,t), \quad j\ge 2,
\end{align*}
as $M=CN$. Then, we  deduce the error bound
\begin{align*}
& \left \vert (S_2-S_2^N)(x_i,t_j) \right \vert
 \le C M^{-1} \sum_{n=1}^j\vert {\cal T}_{S_2;i,n} \vert  \le C M^{-1} + C M^{-1} \sum_{n=2}^j\vert {\cal T}_{S_2;i,n} \vert \nonumber \\
&  \le C M^{-1}+CM^{-1}\sum_{n=2}^j \frac{M^{-1}}{\sqrt{\vr}} E_\gamma(x,t)+CM^{-3/2} \sum_{n=2}^j \frac1{\sqrt{n-1}}+C  M^{-1}  \vr  \sum_{n=2}^j  \frac1{n-1} \nonumber \\
&  \le CM^{-1} +C M^{-3/2} \int_{s=1}^j\frac{ds}{\sqrt{s}}+ C M^{-1}  \vr \int_{s=1}^j\frac{ds}{s} \nonumber \\
&  \le C M^{-1}+ C M^{-1}  \vr \ln M \\
&  \le C M^{-1}.
\end{align*}
Finally, we consider  the error due to the  singular component $ S _1$.  The argument splits into the two cases of $\ve \leq CM^{-1}$ and $\ve \geq CM^{-1}$.  If $ M \ve \geq C$,  from (\ref{S1-bounds}) we obtain the following truncation errors bounds at the first time level $t=t_1$
\[
\vert {\cal T}_{S_1;i,1} \vert  \le  \vr \left \Vert  \frac{\partial^2 S_1}{\partial x^2} \right \Vert_{(x_{i-1},x_{i+1})}
+ a \left \Vert  \frac{\partial S_1}{\partial x} \right \Vert_{(x_{i-1},x_i)}
+  \left \Vert  \frac{\partial S_1}{\partial t} \right \Vert_{(t_0,t_1)} \le C.
\]
At the next time levels $t_n, \, n\ge 2$,  use again (\ref{S1-bounds}) to deduce the truncation error bounds
\begin{align*}
\vert {\cal T}_{S_1;i,j} \vert & \le C \vr (h_i+h_{i+1}) \left \Vert  \frac{\partial^3 S_1}{\partial x^3} \right \Vert_{(x_{i-1},x_{i+1})}
+a h_i \left \Vert  \frac{\partial^2 S_1}{\partial x^2} \right \Vert_{(x_{i-1},x_i)}
+ C k \left \Vert  \frac{\partial^2 S_1}{\partial t^2} \right \Vert_{(t_{j-1},t_j)} \\
& \le C  \left(\frac{N^{-1}}{\vr}+\frac{N^{-1}+M^{-1}}{\sqrt{\vr t_{j-1}}}+\frac{M^{-1}}{t_{j-1}} \right) E_\gamma(x,t)+ CN^{-1}\\
& \le C  \left(\frac{M^{-1}}{\vr}+\frac{M^{-1/2}}{\sqrt{\vr (j-1)}}+\frac1{j-1} \right) E_\gamma(x,t) + CN^{-1}, \quad j\ge 2,
\end{align*}
as $M=CN$. Then,
\begin{align}
& \left \vert (S_1-S_1^N)(x_i,t_j) \right \vert
 \le C M^{-1} \sum_{n=1}^j\vert {\cal T}_{S_1;i,n} \vert
 \le C M^{-1} + C M^{-1} \sum_{n=2}^j\vert {\cal T}_{S_1;i,n} \vert \nonumber \\
&  \le C M^{-1}+C\frac{M^{-1}}{\sqrt{\vr}} \sum_{n=2}^j \frac{M^{-1}}{\sqrt{\vr}} E_\gamma(x,t)+C \frac{M^{-3/2}}{\sqrt{\vr}}  \sum_{n=2}^j \frac1{\sqrt{n-1}}+C M^{-1}  \sum_{n=2}^j  \frac1{n-1} \nonumber \\
&  \le C\frac{M^{-1}}{\sqrt{\vr}} +C \frac{M^{-3/2}}{\sqrt{\vr}} \int_{s=1}^j\frac{ds}{\sqrt{s}}+ C M^{-1} \int_{s=1}^j\frac{ds}{s}, \nonumber \\ 
&  \le  C M^{-1/2} + C M^{-1}\ln M \nonumber \\
&  \le C M^{-1/2}. \label{ErrorS1-Case1}
\end{align}
In the other case of $ M \ve \leq C$,  from~\eqref{psi1-derivatives} we first note the following bounds
\[
\left \vert \frac{\partial \psi ^+_1}{\partial x } \right \vert +  \left \vert \frac{\partial ^2 \psi ^+_1}{\partial x ^2} \right \vert \leq  \frac{C}{\ve} E_\gamma(x,t),
\]
and from \eqref{L-psi+} and \eqref{bounds-psi12}
\[
 \left \vert \frac{\partial \psi ^+_1}{\partial t  } \right \vert  \le \vert L \psi ^+_1 \vert + \vr  \left \vert \frac{\partial ^2 \psi ^+_1}{\partial x ^2} \right \vert + a(t) \left \vert \frac{\partial \psi ^+_1}{\partial x } \right \vert  \leq  C E_\gamma(x,t),
\]
which should be compared to  (\ref{S1-bounds}). Also, $a(0)LS_1 = L \psi ^+_1$. Hence,
\begin{align*}
\vert {\cal T}_{S_1;i,j} \vert  & \le C \vr \left \Vert  \frac{\partial^2 \psi ^+_1}{\partial x^2} \right \Vert_{(x_{i-1},x_{i+1})}
+ C\left \Vert  \frac{\partial \psi ^+_1}{\partial x} \right \Vert_{(x_{i-1},x_i)}
+ C \left \Vert  \frac{\partial \psi ^+_1}{\partial t} \right \Vert_{(t_{j-1},t_j)}  \\
& \le CE_\gamma(x,t).
\end{align*}
We now have
\begin{align}
\left \vert (S_1-S_1^N)(x_i,t_j) \right \vert
& \le C M^{-1} \sum_{n=1}^j\vert {\cal T}_{S_1;i,n} \vert \le C  \sqrt{\vr} \sum_{n=1}^j \frac{M^{-1}}{\sqrt{\vr}}E_\gamma(x,t) \nonumber \\
& \le C \sqrt{\vr} \sum_{n=1}^j \frac{M^{-1}}{\sqrt{\vr}} e^{-\gamma \frac{(x_i-at_j)^2}{4\vr T}} \le C \sqrt{\vr} \nonumber \\
& \le  C M^{-1/2}, \label{ErrorS1-Case2}
\end{align}
where we have used that $\int_{r=-\infty}^\infty \frac1{p}e^{-\frac{r^2}{p}} dr=\sqrt{\pi}$. From~\eqref{ErrorS1-Case1} and~\eqref{ErrorS1-Case2}, we deduce
\[
\left \vert (S_1-S_1^N)(x_i,t_j) \right \vert \le C M^{-1/2}.
\]
Combining all of the bounds above, we deduce the nodal error bound
\begin{eqnarray*}
\Vert    (Y- y)(x_i,t_j) \Vert _{\bar G^{N,M}} \leq  C \vert A_1 \vert N^{-1/2}+C N^{-1} \ln N.
\end{eqnarray*}

Use the arguments in \cite{cd-disc-initialA} to extend this nodal error bound to the global error bound.

\end{proof}

\begin{remark}\label{bdy-int-meet}
If the convective coefficient $a(t)$ only depends on the time variable and the constraint (\ref{1D}) is not imposed on the final time $T$, then the interior layer will interact with the boundary layer (see \cite{asymptotic1, cd-disc-initialA})  in an $O(\sqrt{\ve})$ neighbourhood of the point
$(1,T_*)$, where $d(T_*)=1$. To retain the parameter uniform error bound  (as stated in Theorem \ref{th_a(t)}),
an additional piecewise uniform Shishkin mesh in time should be used either side of $t=T_*$. See  \cite{cd-disc-initialA} for details of the mesh and the associated proof of uniform convergence. Minor modifications to the proof of the error bound are required to deal with the presence of additional terms involving $\psi ^+_i(1,t), i=0,1,2,3,4.$  Example \ref{ex3f} in the numerical section deals with this case of the interior layer and boundary layer interacting.
\end{remark}

In the final theorem, we consider the case of $g_L(0)=\phi (0)$, where the solution  of (\ref{cont-problem-general}) is continuous. In this case the solution of problem~\eqref{cont-problem-general} with $a=a(x,t)$ can be decomposed as in Theorem~\ref{theorem:decomp-general} and
the numerical method  (\ref{discr-probl}) is applied directly to the problem without separating off the singular function $S_0$. The
proof of Theorem~\ref{th_a(t)} is also valid for the following result.

\begin{theorem} \label{th_a(x,t)}  Assume that $a_x(0,0)=0$, $g_L(0)=\phi (0)$ and $ M=O(N)$.
If $Y$ is the solution of (\ref{discr-probl}) and $u$ is the solution of (\ref{cont-problem-general}), then
\begin{eqnarray*}
\Vert   \bar Y- u \Vert _{\bar G} \leq  C \vert A_1 \vert N^{-1/2}+C N^{-1} \ln N.
\end{eqnarray*}
\end{theorem}
\section{Numerical experiments} \label{sec:5}

\begin{figure}[h!] \centering
\resizebox{\linewidth}{!}{ 	\begin{subfigure}[ Computed component  $Y$]{
		\includegraphics[scale=0.5, angle=0]{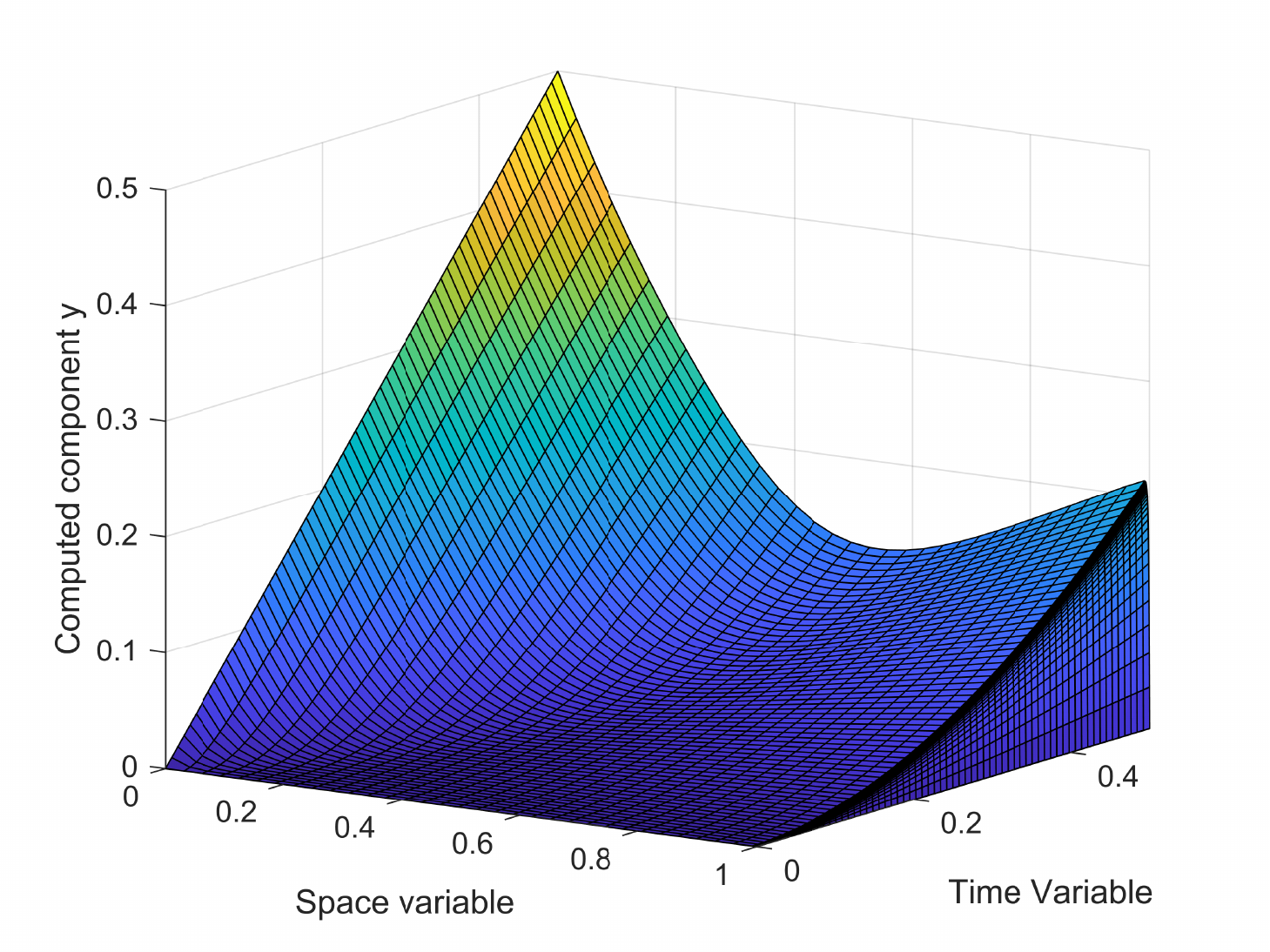}		}    \end{subfigure}
\begin{subfigure}[Numerical approximation $U$]{
		\includegraphics[scale=0.5, angle=0]{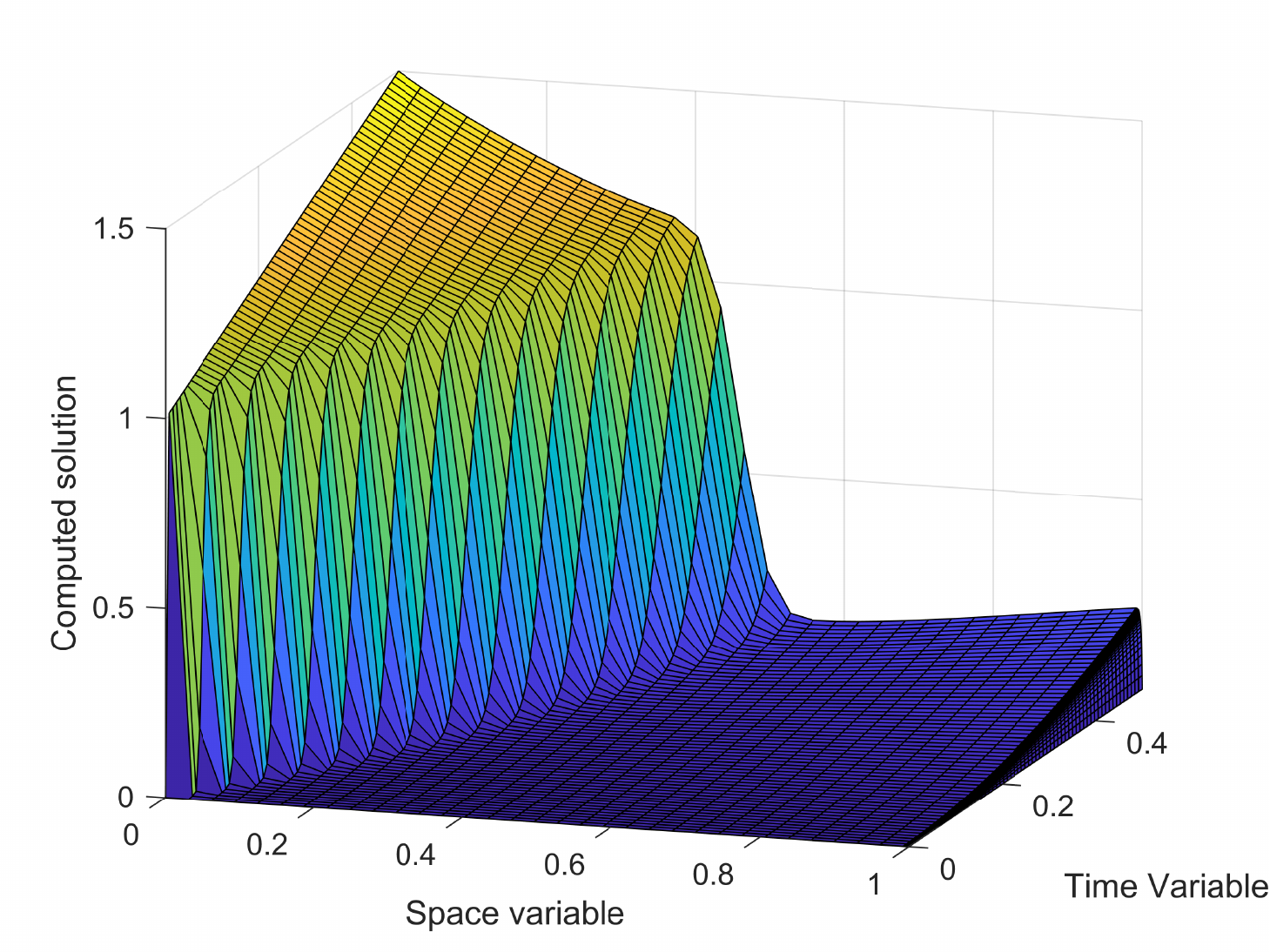}		}	\end{subfigure}}
	\caption{ Example~\ref{ex1f} with $\vr=2^{-10}$: Computed component $Y$ with the scheme~\eqref{discr-probl} for $N=M=64$ and the numerical approximation $U$.}
	\label{fig:Ex1f}
 \end{figure}

The solution of all the  test examples presented below  is unknown and the global orders of convergence are estimated using the two-mesh method~\cite[Chapter 8]{fhmos}. In this particular section, the computed solutions with~\eqref{discr-probl} on the Shishkin meshes $\bar G^{N,M}$ and $\bar G^{2N,2M}$ are denoted, respectively,  by $Y^{N,M}$ and $Y^{2N,2M}$.
Let $\bar Y^{N,M}$  be the bilinear interpolation of the discrete solution $Y^{N,M}$  on the mesh $\bar G^{N,M} $.
Then, compute the maximum two-mesh global differences
$$
D^{N,M}_\ve:= \Vert \bar Y^{N,M}-\bar Y^{2N,2M}\Vert _{\bar G^{N,M} \cup \bar G^{2N,2M}}
$$
and use these values to estimate the orders of global  convergence $ P^{N,M}_\ve$
$$
 P^{N,M}_\ve:=  \log_2\left (\frac{D^{N,M}_\ve}{D^{2N,2M}_\ve} \right).
$$
The uniform  two-mesh global differences $D^{N,M}$ and the uniform orders of global convergence $ P^{N,M}$ are calculated by
$$
D^{N,M}:= \max_{\ve \in S} D^{N,M}_\ve, \quad  P^{N,M}:=  \log_2\left ( \frac{D^{N,M}}{D^{2N,2M}} \right),
$$
where $S=\{2^0,2^{-1},\ldots,2^{-30}\}$.  In all of the tables below  we display the maximum and uniform two-mesh global differences and the corresponding orders of convergence for $N=16,32,\ldots,1024$ and $N=M$.  For the sake of brevity, we display the results in the tables for a smaller representative set of values of $\ve$.
 Note that in the first three examples, the convective coefficient $a(t)$ does not depend on the spatial variable.
\begin{example} \label{ex1f}
We consider the following initial-boundary value problem
\begin{align*}
&u_t-\vr u_{xx}+(1-t^2) u_x=2tx, \quad (x,t)\in (0,1)\times(0,0.5], \\
&u(x,0)=0, \quad x\in (0,1), \\
&u(0,t)=1+t, \ u(1,t)=0, \quad t\in [0,0.5].
\end{align*}
Note that $a'(0) =0$ and $A_1 =1 \ne 0$ in this example.  In Figure~\ref{fig:Ex1f} the computed component $Y$ with the scheme~\eqref{discr-probl} for $\vr=2^{-10}$ and $N=M=64$ is shown. The approximation $U$ to the solution of Example~\ref{ex1f} also appears in that figure; the interior layer emanating from the point $(0,0)$ and the boundary layer in the outflow boundary are observed. The numerical results are in Table~\ref{table:Ex1f} and they  indicate that the numerical method~\eqref{discr-probl} converges uniformly and globally with order $O(N^{-1/2})$ in agreement with Theorem~\ref{th_a(t)}.
\end{example}

\begin{table}[h]
\caption{ Maximum two-mesh global differences and orders of convergence for Example~\ref{ex1f}}
\begin{center}{\tiny \label{table:Ex1f}
\begin{tabular}{|c||c|c|c|c|c|c|c|}
 \hline  & N=M=16 & N=M=32 & N=M=64 & N=M=128 & N=M=256 & N=M=512 & N=M=1024\\
\hline \hline
 $\vr=2^{0}$
&2.593E-03 &1.306E-03 &6.567E-04 &3.285E-04 &1.643E-04 &8.212E-05 &4.106E-05 \\
&0.989&0.992&1.000&1.000&1.000&1.000&
\\ \hline $\vr=2^{-6}$
&3.004E-02 &1.768E-02 &1.014E-02 &5.522E-03 &2.888E-03 &1.467E-03 &7.321E-04 \\
&0.764&0.802&0.877&0.935&0.977&1.003&
\\ \hline $\vr=2^{-12}$
&3.547E-02 &2.356E-02 &1.598E-02 &1.103E-02 &7.581E-03 &5.175E-03 &3.442E-03 \\
&0.590&0.560&0.535&0.541&0.551&0.588&
\\ \hline $\vr=2^{-18}$
&3.551E-02 &2.363E-02 &1.609E-02 &1.118E-02 &7.818E-03 &5.502E-03 &3.877E-03 \\
&0.588&0.554&0.526&0.516&0.507&0.505&
\\ \hline $\vr=2^{-24}$
&3.551E-02 &2.363E-02 &1.609E-02 &1.118E-02 &7.820E-03 &5.505E-03 &3.882E-03 \\
&0.588&0.554&0.525&0.516&0.506&0.504&
\\ \hline $\vr=2^{-30}$
&3.551E-02 &2.363E-02 &1.609E-02 &1.118E-02 &7.820E-03 &5.506E-03 &3.882E-03 \\
&0.588&0.554&0.525&0.516&0.506&0.504&
\\ \hline $D^{N,M}$
&3.551E-02 &2.363E-02 &1.609E-02 &1.118E-02 &7.820E-03 &5.506E-03 &3.882E-03 \\
$P^{N,M}$ &0.588&0.554&0.525&0.516&0.506&0.504&\\ \hline \hline
\end{tabular}}
\end{center}
\end{table}

\begin{example} \label{ex2f}
Consider the example
\begin{align*}
&u_t-\vr u_{xx}+(1-t^2)u_x=2tx, \quad (x,t)\in (0,1)\times(0,0.5], \\
&u(x,0)=x^3, \quad x\in (0,1), \\
&u(0,t)=1+t^2, \ u(1,t)=1, \quad t\in [0,0.5].
\end{align*}
In this example the data problem satisfy that $a'(0) =0$ but $A_1 = 0$. The numerical results obtained for Example~\ref{ex2f} with the numerical method~\eqref{discr-probl}  are given in Table~\ref{table:Ex2f} and they indicate that the method converges with almost first order as stated in Theorem \ref{th_a(t)}.
\end{example}

\begin{table}[h]
\caption{ Maximum two-mesh global differences and orders of convergence for Example~\ref{ex2f}}
\begin{center}{\tiny \label{table:Ex2f}
\begin{tabular}{|c||c|c|c|c|c|c|c|}
 \hline  & N=M=16 & N=M=32 & N=M=64 & N=M=128 & N=M=256 & N=M=512 & N=M=1024\\
\hline \hline
 $\vr=2^{0}$
&6.959E-03 &3.209E-03 &1.541E-03 &7.536E-04 &3.725E-04 &1.852E-04 &9.233E-05 \\
&1.117&1.058&1.032&1.016&1.008&1.004&
\\ \hline $\vr=2^{-6}$
&5.202E-02 &2.956E-02 &1.666E-02 &9.232E-03 &5.067E-03 &2.751E-03 &1.484E-03 \\
&0.816&0.827&0.852&0.865&0.881&0.891&
\\ \hline $\vr=2^{-12}$
&6.938E-02 &3.622E-02 &2.037E-02 &1.125E-02 &6.138E-03 &3.329E-03 &1.792E-03 \\
&0.938&0.830&0.856&0.875&0.883&0.894&
\\ \hline $\vr=2^{-18}$
&6.968E-02 &3.634E-02 &2.044E-02 &1.130E-02 &6.161E-03 &3.342E-03 &1.799E-03 \\
&0.939&0.830&0.855&0.875&0.882&0.894&
\\ \hline $\vr=2^{-24}$
&6.969E-02 &3.634E-02 &2.044E-02 &1.130E-02 &6.162E-03 &3.343E-03 &1.799E-03 \\
&0.939&0.830&0.855&0.875&0.882&0.894&
\\ \hline $\vr=2^{-30}$
&6.969E-02 &3.634E-02 &2.044E-02 &1.130E-02 &6.163E-03 &3.340E-03 &1.802E-03 \\
&0.939&0.830&0.855&0.874&0.884&0.890&
\\ \hline $D^{N,M}$
&6.969E-02 &3.634E-02 &2.044E-02 &1.130E-02 &6.163E-03 &3.343E-03 &1.802E-03 \\
$P^{N,M}$ &0.939&0.830&0.855&0.874&0.883&0.891&\\ \hline \hline
\end{tabular}}
\end{center}
\end{table}

\begin{figure}[h!] \centering
\resizebox{\linewidth}{!}{ 	\begin{subfigure}[ Computed component  $Y$]{
		\includegraphics[scale=0.5, angle=0]{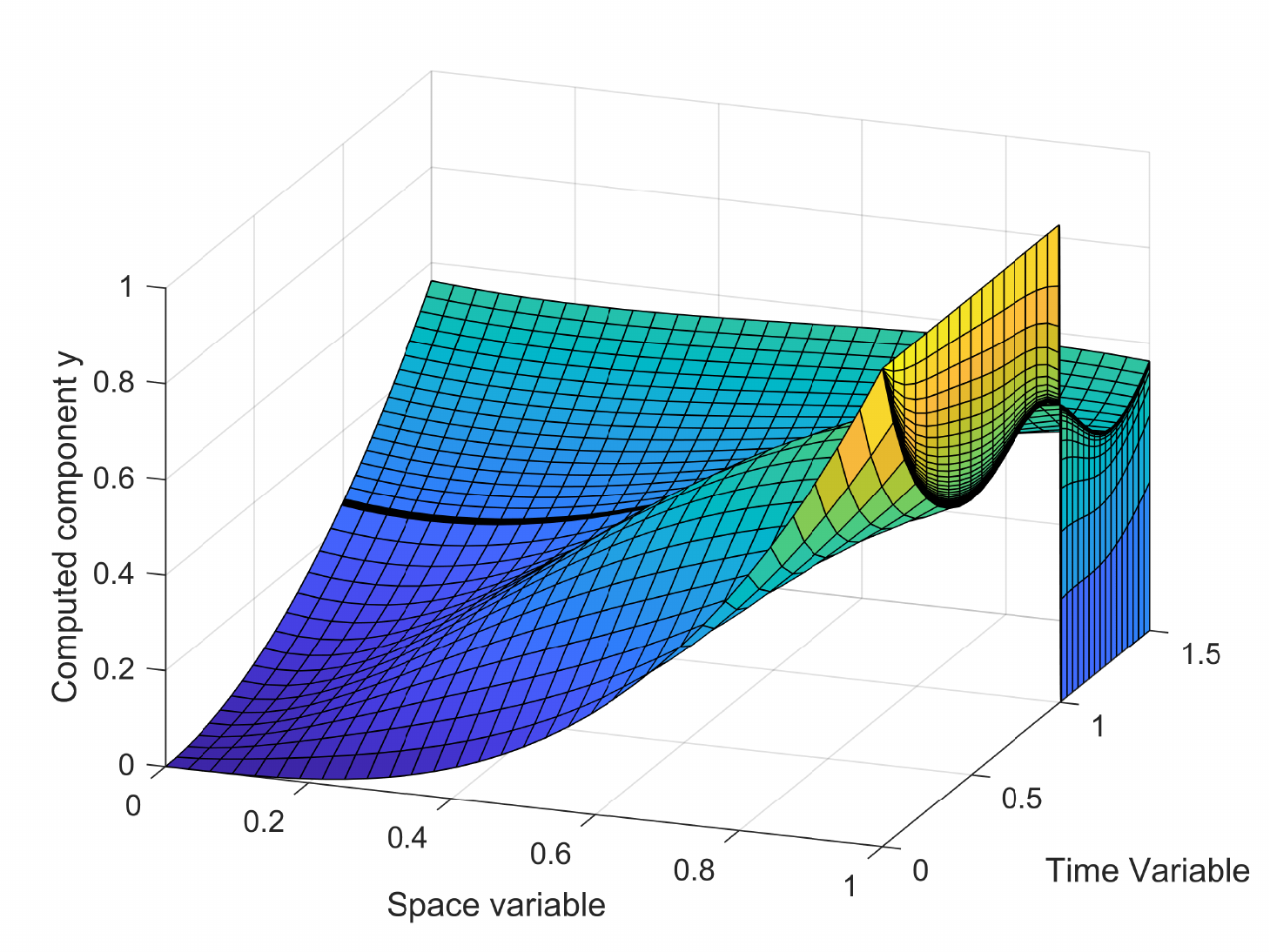}	}    \end{subfigure}
\begin{subfigure}[Numerical approximation $U$]{
		\includegraphics[scale=0.5, angle=0]{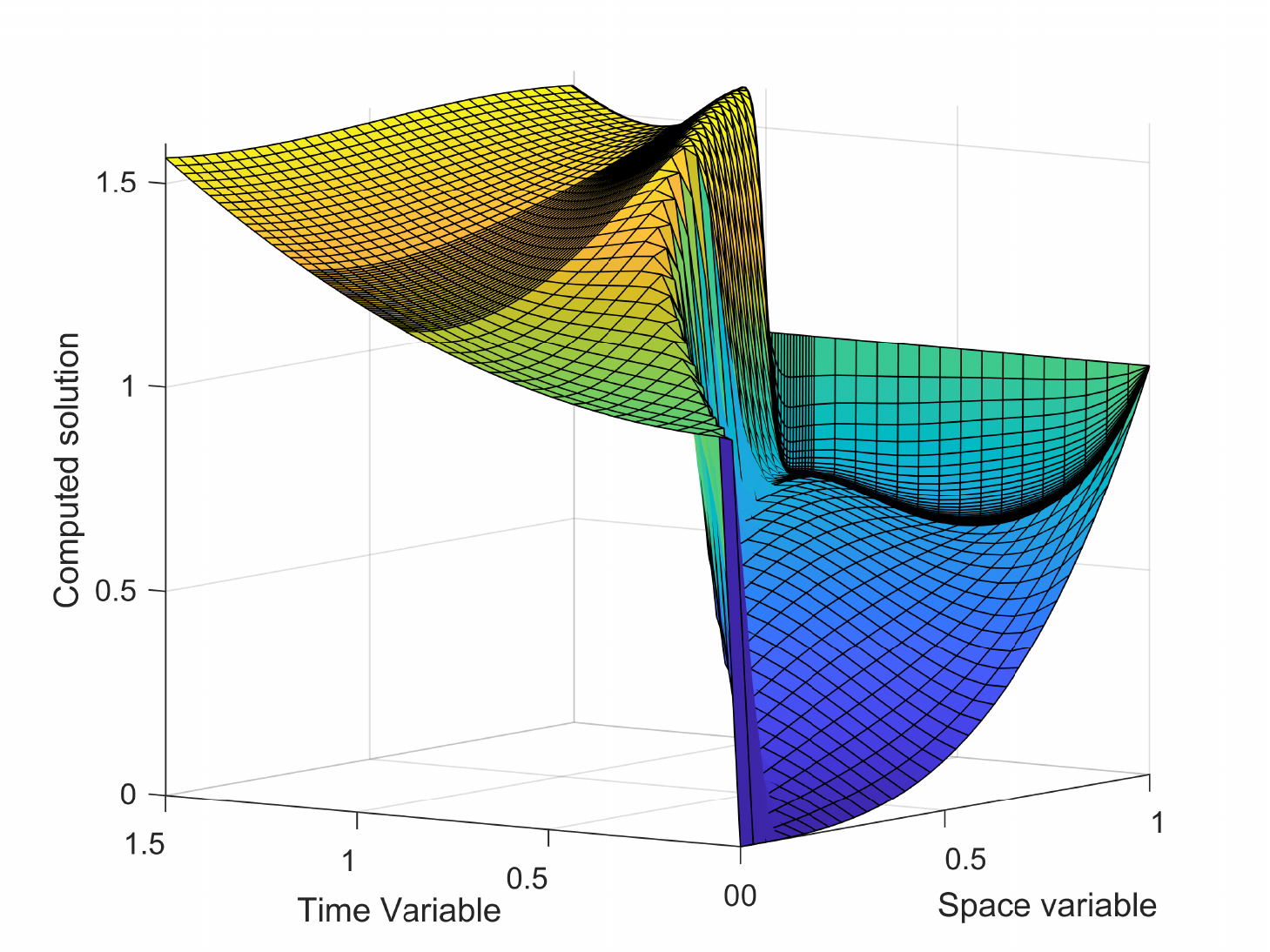}		}	\end{subfigure}}
	\caption{ Example~\ref{ex3f} with $\vr=2^{-10}$: Computed component $Y$ with the scheme~\eqref{discr-probl} for $N=M=64$ and the numerical approximation $U$.}
	\label{fig:Ex3f}
 \end{figure}

\begin{example} \label{ex3f}
Consider the example
\begin{align*}
&u_t-\vr u_{xx}+(1+3t^2-2t)u_x=4x(1-x), \quad (x,t)\in (0,1)\times (0, 1.5], \\
&u(x,0)=x^3, \quad x\in (0,1), \\
&u(0,t)=1+0.25t^2, \ u(1,t)=1, \quad t\in [0,1.5].
\end{align*}
Note that $a'(0) \ne 0, A_1 = 0$ and \eqref{1D} is not satisfied;  then the interior and boundary layers interact with each other. This effect is observed in Figure~\ref{fig:Ex3f} where the approximations to the component $y$ and the solution $u$ are shown.  The numerical results obtained with the scheme~\eqref{discr-probl}
combined with a modification to the mesh in time \cite[(19)]{cd-disc-initialA} {\cbm (see also Remark \ref{bdy-int-meet})}
are given in Table~\ref{table:Ex3f}. These results suggest that the method converges globally and uniformly with almost first order.
\end{example}

\begin{table}[h]
\caption{ Maximum two-mesh global differences and orders of convergence for Example~\ref{ex3f}}
\begin{center}{\tiny \label{table:Ex3f}
\begin{tabular}{|c||c|c|c|c|c|c|c|}
 \hline  & N=M=16 & N=M=32 & N=M=64 & N=M=128 & N=M=256 & N=M=512 & N=M=1024\\
\hline \hline
 $\vr=2^{0}$
&7.807E-02 &3.389E-02 &1.726E-02 &8.507E-03 &4.241E-03 &2.117E-03 &1.057E-03 \\
&1.204&0.973&1.021&1.004&1.002&1.002&
\\ \hline $\vr=2^{-6}$
&1.042E-01 &8.151E-02 &5.318E-02 &2.894E-02 &1.705E-02 &1.018E-02 &5.732E-03 \\
&0.355&0.616&0.878&0.764&0.744&0.828&
\\ \hline $\vr=2^{-12}$
&1.029E-01 &8.090E-02 &5.303E-02 &2.889E-02 &1.702E-02 &1.017E-02 &5.728E-03 \\
&0.347&0.609&0.876&0.763&0.743&0.828&
\\ \hline $\vr=2^{-18}$
&1.025E-01 &8.057E-02 &5.295E-02 &2.887E-02 &1.747E-02 &1.052E-02 &5.982E-03 \\
&0.347&0.606&0.875&0.725&0.732&0.814&
\\ \hline $\vr=2^{-24}$
&1.025E-01 &8.052E-02 &5.293E-02 &2.886E-02 &1.804E-02 &1.094E-02 &6.275E-03 \\
&0.348&0.605&0.875&0.678&0.721&0.802&
\\ \hline $\vr=2^{-30}$
&1.025E-01 &8.052E-02 &5.293E-02 &2.884E-02 &1.813E-02 &1.100E-02 &6.326E-03 \\
&0.348&0.605&0.876&0.670&0.721&0.798&
\\ \hline $D^{N,M}$
&1.062E-01 &8.509E-02 &5.623E-02 &2.897E-02 &1.813E-02 &1.100E-02 &6.326E-03 \\
$P^{N,M}$ &0.319&0.598&0.957&0.676&0.721&0.798&\\ \hline \hline
\end{tabular}}
\end{center}
\end{table}

\begin{example} \label{ex4f}
Consider the example
\begin{align*}
&u_t-\vr u_{xx}+(1+x^2)u_x=4x(1-x), \quad (x,t)\in (0,1)\times(0,0.5], \\
&u(x,0)=0, \quad x\in (0,1), \\
&u(0,t)=u(1,t)=t^2, \  \quad t\in [0,0.5].
\end{align*}
Note that $a_x(0,0) = 0$ and  $A_0=A_1 = 0$.  The numerical results obtained with the scheme~\eqref{discr-probl} are given in Table~\ref{table:Ex4f} and they indicate that the method converges with almost first order as stated in Theorem \ref{th_a(x,t)}.
\end{example}

\begin{table}[h]
\caption{ Maximum two-mesh global differences and orders of convergence for Example~\ref{ex4f}}
\begin{center}{\tiny \label{table:Ex4f}
\begin{tabular}{|c||c|c|c|c|c|c|c|}
 \hline  & N=M=16 & N=M=32 & N=M=64 & N=M=128 & N=M=256 & N=M=512 & N=M=1024\\
\hline \hline
 $\vr=2^{0}$
&2.960E-03 &1.515E-03 &7.615E-04 &3.819E-04 &1.912E-04 &9.567E-05 &4.785E-05 \\
&0.967&0.992&0.996&0.998&0.999&1.000&
\\ \hline $\vr=2^{-6}$
&2.590E-02 &1.643E-02 &9.829E-03 &5.503E-03 &3.026E-03 &1.634E-03 &8.757E-04 \\
&0.656&0.742&0.837&0.863&0.889&0.900&
\\ \hline $\vr=2^{-12}$
&3.115E-02 &2.004E-02 &1.201E-02 &6.698E-03 &3.677E-03 &1.992E-03 &1.070E-03 \\
&0.636&0.739&0.843&0.865&0.884&0.897&
\\ \hline $\vr=2^{-18}$
&3.126E-02 &2.011E-02 &1.205E-02 &6.718E-03 &3.689E-03 &1.999E-03 &1.073E-03 \\
&0.637&0.739&0.843&0.865&0.884&0.897&
\\ \hline $\vr=2^{-24}$
&3.126E-02 &2.011E-02 &1.205E-02 &6.719E-03 &3.689E-03 &1.999E-03 &1.073E-03 \\
&0.637&0.739&0.843&0.865&0.884&0.897&
\\ \hline $\vr=2^{-30}$
&3.126E-02 &2.011E-02 &1.205E-02 &6.718E-03 &3.690E-03 &1.998E-03 &1.075E-03 \\
&0.637&0.739&0.843&0.865&0.885&0.894&
\\ \hline $D^{N,M}$
&3.126E-02 &2.011E-02 &1.205E-02 &6.719E-03 &3.690E-03 &1.999E-03 &1.075E-03 \\
$P^{N,M}$ &0.637&0.739&0.843&0.865&0.884&0.895&\\ \hline \hline
\end{tabular}}
\end{center}
\end{table}

\begin{example} \label{ex5f}
Consider the example
\begin{align*}
&u_t-\vr u_{xx}+(1+x)u_x=4x(1-x), \quad (x,t)\in (0,1)\times(0,0.5], \\
&u(x,0)=0, \quad x\in (0,1), \\
&u(0,t)=t;\ u(1,t)=t^2, \  \quad t\in [0,0.5].
\end{align*}
Note that $a_x(0,0) \neq 0$ and  $A_0=0, A_1 \neq 0$. The numerical results obtained with the scheme~\eqref{discr-probl} are given in Table~\ref{table:Ex5f}. They suggest that the method converges globally and uniformly with order $O(N^{-1/2})$, but the theoretical justification of these results remains open,  as the proof in Theorem~\ref{th_a(x,t)} requires $a_x(0,0)=0$.
\end{example}

\begin{table}[h]
\caption{ Maximum two-mesh global differences and orders of convergence for Example~\ref{ex5f}}
\begin{center}{\tiny \label{table:Ex5f}
\begin{tabular}{|c||c|c|c|c|c|c|c|}
 \hline  & N=M=16 & N=M=32 & N=M=64 & N=M=128 & N=M=256 & N=M=512 & N=M=1024\\
\hline \hline
 $\vr=2^{0}$
 &1.426E-03 &8.292E-04 &4.557E-04 &2.388E-04 &1.222E-04 &6.173E-05 &3.099E-05 \\
&0.782&0.863&0.932&0.967&0.985&0.994&
\\ \hline $\vr=2^{-6}$
&2.563E-02 &1.562E-02 &9.145E-03 &5.113E-03 &2.896E-03 &1.600E-03 &8.695E-04 \\
&0.714&0.772&0.839&0.820&0.856&0.879&
\\ \hline $\vr=2^{-12}$
&3.099E-02 &2.117E-02 &1.437E-02 &9.917E-03 &6.914E-03 &4.794E-03 &3.258E-03 \\
&0.550&0.559&0.535&0.520&0.528&0.557&
\\ \hline $\vr=2^{-18}$
&3.108E-02 &2.128E-02 &1.449E-02 &1.007E-02 &7.115E-03 &5.061E-03 &3.603E-03 \\
&0.547&0.554&0.525&0.501&0.492&0.490&
\\ \hline $\vr=2^{-24}$
&3.108E-02 &2.128E-02 &1.449E-02 &1.007E-02 &7.119E-03 &5.065E-03 &3.609E-03 \\
&0.547&0.554&0.525&0.501&0.491&0.489&
\\ \hline $\vr=2^{-30}$
&3.108E-02 &2.128E-02 &1.449E-02 &1.007E-02 &7.119E-03 &5.065E-03 &3.609E-03 \\
&0.547&0.554&0.525&0.501&0.491&0.489&
\\ \hline $D^{N,M}$
&3.108E-02 &2.128E-02 &1.449E-02 &1.007E-02 &7.119E-03 &5.065E-03 &3.609E-03 \\
$P^{N,M}$ &0.547&0.554&0.525&0.501&0.491&0.489&\\ \hline \hline
\end{tabular}}
\end{center}
\end{table}

 \begin{remark} \label{remark:Numer-H}
 In the numerical experiments performed, we have considered~\eqref{psi0+H}  when evaluating $\psi_0^+$ and $\theta:=d(t_j)x_i/(t_j \vr)$ is a large number with $(x_i,t_j)\in G^{N,M}$ to prevent overflow problems. If $\theta \ge 300$, the value of the  Mill's ratio $H$ in Examples~\ref{ex1f}, \ref{ex2f} and \ref{ex3f} has been computed  by using that
\[
H(r) \sim \frac1{r\sqrt{\pi}} \left(1+\sum_{m=1}^\infty (-1)^m \frac{1.3\ldots(2m-1)}{(2r^2)^m}\right), \text{ as } r \to \infty.
\]
This series has been approximated by the $n$-th partial sum and the maximum two-mesh global differences in Tables~\ref{table:Ex1f}, ~\ref{table:Ex2f} and ~\ref{table:Ex3f} have been obtained with $n=5.$ Similar results have been obtained if larger values of $n$ are considered.
\end{remark}

\section{ Appendix: Bounds on the partial derivatives of the functions $\psi ^\pm_i(x,t)$.} \label{sec:appendix1}

 The functions $\psi ^\pm_i(x,t)$ are defined by means of the iterated integrals of the complementary error function. Define
\[
\erfc_{-1}(x):=\frac{2}{\sqrt{\pi}}e^{-x^2}, \quad \erfc_n(x):=\int_{s=x}^\infty\erfc_{n-1}(s) \, ds, \quad n \geq 0.
\]
 Note that $\erfc _0(x) = \erfc (x)$
\footnote{ The first three iterated integrals of the complementary error function are
\begin{align*}
\erfc _1(x) &= \frac{e^{-x^2}}{\sqrt{\pi}} - x \erfc (x) = e^{-x^2} \left( \frac{1}{\sqrt{\pi}}  - x e^{x^2}\erfc (x) \right),\\
\erfc _2 (x) &=\frac{1}{4} \left( (1+2x^2) \erfc (x) - \frac{2xe^{-x^2}}{\sqrt{\pi}}  \right), \\
 \erfc _3 (x) &= \frac{1}{6} \left( \frac{(1+x^2)e^{-x^2}}{\sqrt{\pi}} - \frac{ (3x+2x^3) }{2}\erfc (x) \right).
\end{align*}
}
 and
\[
\erfc _n (x)=\frac2{\sqrt{\pi}} \int_{s=x}^\infty \frac{(s-x)^n}{n!} e^{-s^2}ds.
\]
 In addition, we have the following identities~\cite{NIST} 
\begin{subequations}\label{identities}
\begin{align}
n\erfc _n (x) +x \erfc _{n-1} (x) &= \frac{1}{2} \erfc _{n-2} (x);\quad n \geq 1; \label{identities-a} \\
(-1)^n\erfc _n (x) +\erfc _n (-x) &=  \frac{i^{-n}}{2^{n-1} n!} H_n(ix) ; \quad n \geq 0,
\end{align}
\end{subequations}
 where $i^2=-1$ and $H_n$ is the Hermite polynomial of degree $n$.
Recall the definitions in  \eqref{def-basic} and note that
\begin{subequations}
\begin{align}
\psi ^\pm _1(x,t) &= (x\pm d(t))\psi _0^\pm (x,t) +2\ve t \psi _{-1}^\pm (x,t), \   \label{initial-recurrance} \\
2t\frac{\partial \psi ^- _1} {\partial t}&=  \psi ^- _1 - \bigl( 2ta(d(t),t)+ (x-d(t))\bigr) \psi ^- _0,\\
2t\frac{\partial \psi _1^+} {\partial t}&= \left(1 + \frac{2xp(t)}{\ve t} \right)\psi _1^+ + \bigl( 2ta(d(t),t)- (x+d(t))\bigr) \psi _0^+,
\label{tdpsi+dt}
\end{align}
where
\begin{equation}\label{def-p}
p(t):= ta(d(t),t)-d(t)  =  \int _{s=0}^t a(d(t),t) - a(d(s),s)) \ ds.
\end{equation}
\end{subequations}
Observe that $p(t) \equiv 0$, when $a(x,t)=a$ is a constant.

 Some recurrence relations are given below which are useful when bounding the derivatives of the functions $\psi ^\pm_n(x,t)$.
\begin{subequations}\label{deriv-relations-general}
For all $n\geq 1$
\begin{equation}
 \frac{\partial \psi ^-_n}{\partial x}  =  n\psi ^- _{n-1}, \quad
 \frac{\partial \psi ^+_n}{\partial x}  =  n\psi ^+_{n-1} + \frac{d(t)}{\ve t} \psi _n^+,
\end{equation}
and for all  $n\geq 2$ (using (\ref{identities}a)) we have
\begin{align}
\psi ^\pm _n(x,t) &=  (x\pm d(t))\psi _{n-1}^\pm (x,t) +2(n-1)\ve t \psi _{n-2}^\pm (x,t), \label{recurrance}
\\
 \frac{\partial \psi ^-_n}{\partial t} &= 
\ve  n(n-1)\psi ^- _{n-2}  -a(d(t),t)  n\psi ^-_{n-1}, \nonumber
\\
 \frac{\partial \psi ^+ _n}{\partial t}
&= \ve  n(n-1)\psi ^+_{n-2} + \left(2\frac{d(t)}{t} -a(d(t),t) \right) n\psi ^+_{n-1} \nonumber
\\
&-  \frac{p(t)}{\ve t^2} (d(t)\psi ^+_{n} - \psi ^+_{n+1}). \nonumber
\end{align}
\end{subequations}

In the case of constant coefficients one has $L\psi ^-_n = L\psi ^+_n =0$, but for variable $a(x,t)$, by using (\ref{identities}a) we have that  for all $n\geq 0$
\begin{subequations}
\begin{align}
L \psi ^- _n  &=  (a(d(t),t)-a(x,t)) \frac{\partial \psi ^- _n}{\partial x}, \label{L-psi}
\\
 L \psi ^+_n &=  (a(d(t),t)-a(x,t)) \frac{\partial \psi ^+_n}{\partial x} + p(t)\frac{ \psi ^+_{n+1}}{\ve t^2}. \label{L-psi+}
\end{align}
\end{subequations}



Using the inequality $\erfc(z) \leq C e^{-z^2} \leq C e^{\gamma ^2/4}e^{-\gamma z}, \forall z \geq 0$
it follows (see \cite{cd-disc-initialA} and \cite{cd-disc-initialB}) that
\begin{subequations}\label{bounds-singular-functions}
\begin{align}
 & \left \vert \frac{\partial ^j }{\partial t ^j }   \psi ^-_0(x,t)\right \vert, \left \vert \frac{\partial ^j }{\partial t ^j }  E(x,t)\right \vert  \leq  C
 \left( \frac{1}{t} +\frac{1}{\sqrt {\ve t}}\right)^j E_\gamma (x,t);\quad j=1,2, \label{bound-A}
 \\
 & \left \vert \psi ^-_0(x,t)\right \vert \leq C \quad \hbox{and} \quad \vert \psi ^- _0(x,t) \vert \leq C E(x,t), \quad \hbox{if} \quad x \geq d(t), \\
 & \left \vert \frac{\partial ^i }{\partial x ^i } \psi ^-_0(x,t)\right \vert, \left \vert \frac{\partial ^i }{\partial x^i }  E(x,t)\right \vert    \leq C\left(\frac{1}{\sqrt {\ve t}}\right) ^{i} E_\gamma (x,t), \quad  1\leq i \leq 4.\label{bound-B}\end{align}
\end{subequations}

The following remark is used  to prove  bounds on the derivatives of the singular function $\psi^+_0$  and  to compute the numerical results  presented in Section \S \ref{sec:5} (see Remark~\ref{remark:Numer-H}).
\begin{remark} \label{remark:H}
The function $\psi^+_0$ can be written as
\begin{equation} \label{psi0+H}
\psi^+_0(x,t)=\frac1{2} e^{\frac{d(t)x}{t\vr}} \erfc \left( \frac{x+d(t)}{2 \sqrt{\vr t}} \right)= \frac1{2} E(x,t) \, H\left( \frac{x+d(t)}{2\sqrt{\vr t}}\right),
\end{equation}
where $H$ is the Mill's ratio and it is defined by
$
H(x):=e^{x^2}\erfc(x).
$
From~\cite{lether}, we have the  inequality
\begin{equation}\label{sharp}
 \frac{1}{\frac{\pi -1}{\sqrt{\pi}} x+\sqrt{1+\frac{x^2}{\pi}}} \leq H(x) \leq \frac{1}{\frac{2}{\sqrt{\pi}} x+\sqrt{1+\frac{(\pi -2)^2x^2}{\pi}}}.
\end{equation}
Hence, for all $x >0$
\[
\left (1-\sqrt{\pi}  xH (x) \right) \leq  \min \left \{ \frac{1}{2x^2}, \frac{1}{\sqrt{\pi}x} \right \},
\]
and
\[
e^{\frac{d(t)x}{t\ve}} \erfc _1 \left( \frac{x+d(t)}{2 \sqrt{\vr t}} \right) = \frac{E(x,t)}{\sqrt{\pi}} \left( 1 -\sqrt{\pi}   \left( \frac{x+d(t)}{2 \sqrt{\ve t}}\right) H\left( \frac{x+d(t)}{2 \sqrt{\ve t}}\right) \right).
\]
Hence,
\begin{equation}\label{boundG}
\frac{1}{\sqrt{\ve t}} \left( e^{\frac{d(t)x}{t\ve}} \erfc _1 \left( \frac{x+d(t)}{2 \sqrt{\vr t}} \right)   \right)
 \leq C \frac{E(x,t)}{x+d(t)}  \min \left \{ 1, \frac{\sqrt{\ve t}}{x+d(t)} \right \}.
\end{equation}
\end{remark}

 In the next lemma bounds on the derivatives of the function $ \psi^+_0$ are deduced.
\begin{lemma} \label{lemma:BoundsPsi0+}
For the singular function $ \psi^+_0$,
we have the following bounds
\begin{subequations}\label{bounds-star-functions}
\begin{align}
 \vert \psi ^+_0 (x,t) \vert & \leq C \min \left \{ 1, \frac{\sqrt{\ve t}}{x+d(t)} \right \}E(x,t),
 \label{bounds-star-functions-a}
 \\
  \left \vert \frac{\partial  }{\partial t }  \psi ^+ _0(x,t) \right \vert   & \leq  \frac{C}{t} E_\gamma (x,t),
 \label{bounds-star-functions-b}
 \\
  \left \vert \frac{\partial ^2 }{\partial t ^2}  \psi ^+ _0(x,t)\right \vert  & \leq
\frac{C}{t^2} \left(1+\sqrt{\frac{t}{\ve}} \right) E_\gamma (x,t),
\label{bounds-star-functions-c}
\\
\left \vert \frac{\partial }{\partial x }  \psi ^+ _0(x,t)\right \vert   & \leq  \frac{C}{x+d(t)}E_\gamma(x,t) ,
 \label{bounds-star-functions-d}
 \\
 \left \vert \frac{\partial ^2}{\partial x ^2}  \psi ^+ _0(x,t)\right \vert   & \leq  \frac{C}{\ve t}E_\gamma(x,t),
 \label{bounds-star-functions-e}
 \\
\left \vert \frac{\partial ^3}{\partial x ^3}  \psi ^+ _0(x,t)\right \vert   & \leq  \frac{C}{\ve ^2t} \left(1+ \sqrt{\frac{\ve}{t}} \right)E_\gamma(x,t). \label{bounds-star-functions-f}
  \end{align}
\end{subequations}
\end{lemma}
\begin{proof}

Note first that
\begin{align*}
p(t)= ta(d(t),t)-d(t)
=  \int _{s=0}^t \int _{r=d(s)}^{d(t)}   a_x(r,t)   dr\,  ds +\int _{s=0}^t \int _{r=s}^t a_t(d(s),r)  dr\,  ds  .
\end{align*}
Hence,
$
\vert p(t)\vert \leq Ct^2.
$\footnote{ In the particular case of $\nabla a(0,0) =(0,0)$, one has $\vert p(t) \vert \leq Ct^3$.  }
To prove~\eqref{bounds-star-functions-a},  we use that $\vert H(r) \vert \leq C$ and  $rH(r) \leq C$ for all $r \geq 0$, where $H$ is defined in Remark \ref{remark:H}.   Then,
\[
(x+d(t))\psi ^+_0(x,t) \leq C\sqrt{\ve t} E(x,t)  \quad \hbox{and} \quad \psi ^+_0 (x,t)  \leq CE(x,t).
\]
  Using \eqref{boundG}, (\ref{bounds-star-functions-a})  and the identity (\ref{recurrance})
 we easily establish the following bounds
\begin{equation}\label{bounds-psi12}
\,  \vert \psi ^+ _i(x,t) \vert \leq  C (\sqrt{\ve t} )^{i} \min \left\{ 1, \sqrt{ \frac{\vr}{ t}} \right\} E(x,t), \quad i=1,2.
\end{equation}
To prove~\eqref{bounds-star-functions-b}, observe that
\begin{align*}
& (x+d(t))\frac{\partial \psi ^+_0(x,t)}{\partial t} = \frac{xp(t)}{t}\left(\frac{(x+d(t))}{t\ve }\psi ^+_0(x,t) -\frac{E(x,t)}{\sqrt{\ve \pi t} }  \right)
\nonumber \\
& \hspace{1cm} +\left(4xp(t) +(x+d(t))((x+d(t)-2ta(t)) \right) \frac{E(x,t)}{4t\sqrt{\ve \pi t} } \nonumber
\\
 &  \hspace{1cm} = \frac{xp(t)}{\ve t^2} \psi ^+_1(x,t)  +
(x-d(t))^2 \frac{E(x,t)}{4t\sqrt{\ve \pi t}}  +a(d(t),t) (x-d(t))  \frac{E(x,t)}{2\sqrt{\ve \pi t} }, \nonumber
\end{align*}
 and use (\ref{boundG}) and $\vert p(t) \vert \leq Ct^2$.  
 Next, we prove \eqref{bounds-star-functions-c}.
We have that
\begin{align*}
(x+d(t))\frac{\partial ^2 \psi ^+_0(x,t)}{\partial t^2} & = \frac{\partial}{\partial t} \Bigl(  \frac{xp(t)}{\ve t^2} \psi ^+_1 (x,t) +
(x+d(t)  +2p(t))  \frac{ (x-d(t))E(x,t)}{4t\sqrt{\ve \pi t} }\Bigr) \\
& - a(d(t),t) \frac{\partial \psi ^+_0(x,t)}{\partial t}
\end{align*}
\begin{align*}
\hbox{and} &\left \vert \frac{\partial }{\partial t}\left( \frac{p(t)}{t^2} \right) \right \vert \leq C +C \frac{ \Vert \nabla a (0,0) \Vert }{t}  \le \frac{C}{t},
\\
& \left \vert \frac{(x-d(t))}{2\sqrt{\ve t}} \frac{\partial E(x,t)}{\partial t} \right \vert \leq \frac{C}{t} \left(1+\sqrt{\frac{t}{\ve}}\right)E_\gamma(x,t) ,
\\
& \left \vert \frac{\psi ^+_1(x,t)}{\ve t} \right \vert \leq \frac{C}{t} \min \left \{ 1, \sqrt{\frac{t}{\ve} } \right \} E_\gamma(x,t),
\\
& \frac{\partial \psi _1^+(x,t)} {\partial t}= \left(\frac{1}{2t} + \frac{xp(t)}{\ve t^2} \right)\psi _1^+(x,t) + \frac{2p(t)- (x-d(t))}{2t}\psi _0^+(x,t).
\end{align*}
Collecting all of these bounds yields \eqref{bounds-star-functions-c}.
 From
\[
(x+d(t))\frac{\partial }{\partial x} \psi ^+_0 (x,t)= d(t) \left( \frac{x+d(t)}{t\ve }\psi ^+_0(x,t)  -\frac{E(x,t)}{\sqrt{\ve \pi t} } \right) -\frac{(x-d(t))E(x,t)}{2\sqrt{\ve \pi t} },
\]
 and~\eqref{boundG}, we have ~\eqref{bounds-star-functions-d}.
Note also that for $i=2,3$
\[
\frac{\partial ^i}{\partial x ^i}\psi ^+_0(x,t) = \frac{d(t)}{t\ve} \frac{\partial ^{i-1}}{\partial x ^{i-1}} \psi ^+_0(x,t) - \frac{1}{2\sqrt{\ve \pi t}}\frac{\partial ^{i-1}}{\partial x ^{i-1}}E(x,t),
\]
from which~\eqref{bounds-star-functions-e} and~\eqref{bounds-star-functions-f}  follows.
\end{proof}

\end{document}